\documentclass[envcountsect, envcountsame, twocolumn, amsthm]{autart}
\usepackage[utf8x]{inputenc}
\usepackage[english]{babel}
\usepackage{amsmath}
\usepackage{amssymb}
\usepackage{xcolor}
\usepackage[colorlinks=false,linkbordercolor=red,citebordercolor=blue]{hyperref}
\usepackage{tabularx}
\usepackage{comment}
\usepackage{graphicx}
\usepackage{wrapfig}
\usepackage{epsfig}
\usepackage{caption}
\usepackage{subcaption}
\usepackage{enumitem}

\usepackage{adjustbox}

\newtheorem{Theorem}{Theorem}[section]

\newtheorem{Definition}[Theorem]{Definition}
\newtheorem{Lemma}[Theorem]{Lemma}

\newtheorem{Remark}[Theorem]{Remark}

\newcommand{\setdef}[2]{\left\{\, #1\, \left|\, \vphantom{#1} #2 \right.\right\}}

\newcommand{\ddt}{\tfrac{\text{\normalfont d}}{\text{\normalfont d}t}}
\newcommand{\dt}{\text{\normalfont d}t}

\newcommand{\R}{\mathbb{R}}

\newcommand{\N}{\mathbb{N}}

\newcommand{\cC}{\mathcal{C}}
\newcommand{\cD}{\mathcal{D}}
\newcommand{\cF}{\mathcal{F}}

\newcommand{\cL}{\mathcal{L}}
\newcommand{\cN}{\mathcal{N}}

\newcommand{\cT}{\mathcal{T}}

\newcommand{\cW}{\mathcal{W}}

\newcommand{\vp}{\varphi}
\newcommand{\ve}{\varepsilon}
\newcommand{\rp}{\mathbb{R}_{\geq 0}}

\newcommand{\oT}{\textbf{T}}

{\left(\begin{smallmatrix}}
{\end{smallmatrix}\right)}

{\left[\begin{smallmatrix}}
{\end{smallmatrix}\right]}

\DeclareMathOperator{\Gl}{\mathbf{Gl}}

\DeclareMathOperator{\esssup}{\rm ess\,sup}

\usepackage{tikz}
\usetikzlibrary{backgrounds}
\usetikzlibrary{intersections}
\usetikzlibrary{shapes, arrows, math, calc, shapes.multipart}
\tikzstyle{block} = [draw, fill=white, rectangle, minimum height=3em, minimum width=6em]
\tikzstyle{output} = [coordinate]
\tikzstyle{input} = [coordinate]
\usetikzlibrary{decorations.markings,arrows.meta}
\usetikzlibrary{decorations.pathmorphing,shapes,arrows}

\usepackage{pagecolor}
\usetikzlibrary{decorations, decorations.text,backgrounds}
\tikzstyle{int}=[draw, fill=blue!20, minimum size=2em]
\tikzstyle{init} = [pin edge={to-,thin,black}]
\usepackage{pgfplots}
\usepackage{tkz-euclide}

\newcommand\blfootnote[1]{%
  \begingroup
  \renewcommand\thefootnote{}\footnote{#1}%
  \addtocounter{footnote}{-1}%
  \endgroup
}

\pgfplotsset{compat=1.18}

\begin{document}

\begin{frontmatter}
\title{Exact output tracking in prescribed finite time via funnel control\thanksref{footnoteinfo}} 

\thanks[footnoteinfo]{
This work was supported by the German Research Foundation (DFG, Deutsche Forschungsgemeinschaft); 
Project-IDs: 362536361, and 471539468.
}

\author[UPB,TUI]{Lukas Lanza}\ead{{lukas.lanza@tu-ilmenau.de}} 

\address[UPB]{Institute of Mathematics, University of Paderborn, Warburger Str. 100, 33098 Paderborn, Germany}  
\address[TUI]{
Optimization-based Control Group, Institute of Mathematics, Technische Universit\"at Ilmenau, Germany
}


\begin{abstract}
Output reference tracking of unknown nonlinear systems is considered.
The control objective is exact tracking in predefined finite time, while in the transient phase the tracking error evolves within a prescribed boundary.
To achieve this, a novel high-gain feedback controller is developed that is similar to, but extends, existing high-gain feedback controllers.
Feasibility and functioning of the proposed controller is proven rigorously.
Examples for the particular control objective under consideration are, for instance, linking up two train sections, or docking of spaceships.
\end{abstract}

\end{frontmatter}

\section{Introduction}
The control objective, to bring the output of a system to a certain exact value within predefined finite time has various applications: in modern robot-based industry (placement of components), in public transportation (connection of two train sections), autonomous driving (docking at the charging station), and in astronautics (rendezvous of spacecraft), to name but a few.
While \emph{asymptotic exact tracking} has been studied for some time, there are few results on \emph{exact tracking in finite time}.
In~\cite{IlchRyan02b}, referring to the results in~\cite{ByrnWill84}, it is shown that the proposed funnel controller achieves global asymptotic stabilization for a class of linear multi-input multi-output (MIMO) systems of relative degree one.
A generalization to a class of nonlinear relative degree one MIMO systems is proposed in~\cite{RyanSang09}.
In~\cite{NuneHsu09} an extended sliding mode controller is proposed, which achieves asymptotic tracking of linear single-input single-output (SISO) systems. 
This controller is extended to linear MIMO systems in~\cite{NunePeix14}.
In~\cite{Davi13} backstepping is combined with feedback linearization techniques and higher order sliding modes to design a controller, which achieves exponential tracking.
In~\cite{OlivPeix13} a high-gain based sliding mode controller is introduced, where the peaking related to the high-gain observer is avoided by introducing a dwell-time activation scheme. 
This controller achieves asymptotic tracking for a class of nonlinear SISO systems of arbitrary relative degree, where the reference signal is generated by a reference model.
To the price of a discontinuous control, asymptotic tracking for nonlinear MIMO systems is achieved in~\cite{VergDima19,VergDima21}.
In~\cite{LeeTren19} a funnel controller is proposed, which achieves asymptotic tracking for a class of nonlinear relative degree one MIMO systems.
This result is extended in the recent work~\cite{BergIlch21}, where it is shown that the proposed controller achieves asymptotic tracking of nonlinear MIMO systems with arbitrary relative degree whereas the tracking error has prescribed transient behaviour.
We now turn from asymptotic tracking towards exact tracking in finite time.
In~\cite{Andr08} a recursive observer structure as well as an extension of the homogeneous approximation technique is introduced to achieve global asymptotic as well as finite time stabilization for higher order chain of integrator systems.
For control affine systems with given relative degree, in~\cite{Leva14} a homogeneous higher order sliding mode controller is designed, which achieves stabilization in finite time.
In~\cite{BasiPana16} sliding mode control concepts and results from~\cite{Leva03,BhatBern05} are used to establish control schemes, which achieve finite time stabilization for linear SISO \& MIMO systems.
In~\cite{EstrFrid16} backstepping and higher order sliding mode control are combined to construct a controller, which achieves exact output tracking in finite time for nonlinear MIMO systems in nonlinear block controllable form.
Similar to the prescribed performance controller in~\cite{DimaBech20}, this controller suffers from the proper initialization problem, where it is not clear, how large to choose the involved parameters.
The controller in~\cite{EstrFrid16}, along with limiting conditions on the system class, presumes knowledge of the system's functions and explicitly involves inverses of some.
In~\cite{YangDing18} a controller is introduced, which achieves exact tracking in finite time for a class of nonlinear SISO systems satisfying a certain homogeneity assumption. 
This controller relies on estimations of the external disturbances, where the problem of proper initialization is avoided by assuming explicit knowledge of the bounds of the disturbances and the reference.
The controller explicitly involves (parts of) the system's right hand side and is of relatively high complexity.
Utilizing the implicit Lyapunov function approach, in~\cite{MercUrib21} a state feedback-integral controller is designed, which is capable to stabilize homogeneous systems (negative and positive) in fixed finite time, where the final time can be estimated involving the initial state and a corresponding Lyapunov function.
In most of the control schemes for exact tracking in finite time discussed above, the final time
cannot be prescribed; only the existence of such a finite time is ensured.
Contrary, in~\cite{RodrSanc17} a controller is introduced which solves a \textit{predefined-time exact tracking problem} for the class of fully actuated mechanical (relative degree two) systems.
The controller relies on a backstepping procedure and consists of a \textit{predefined-time stabilization function}, and involves the system's equations explicitly.
In~\cite{Tren19} a funnel controller is introduced, which achieves asymptotic tracking as well as convergence to zero of the tracking error in finite time for a class of relative degree one SISO systems.
Regularization in prescribed finite time is achieved in~\cite{song2017time} for nonlinear systems in strict feedback form, invoking a strictly increasing scaling function for the state, which grows unbounded when approaching the final time.
In the recent work~\cite{EspiPerr22}, linear time invariant systems with delayed input are under consideration. Representing the delay system in a PDE-ODE cascade, under the usage of backstepping techniques and integral transformations a controller is designed, which stabilizes the system within predefined finite time.
\\
Circumventing some drawbacks mentioned above, we propose a controller, which achieves \textit{exact tracking in predefined finite time}.
Since the controller is of funnel type it inherits the advantages of robustness with respect to noise, and that the tracking error evolves within prescribed bounds. 
Moreover, the controller is model-free in the sense that no knowledge of the system's parameters is assumed; 
only knowledge of the order~$r \in \N$ of the differential equation and the common dimension~$m \in \N$ of the input and output is required, and the system's right-hand side has to satisfy a high-gain property.
The system class under consideration is the same as in~\cite{BergIlch21},
and encompasses the systems under consideration in \cite{Andr08,Leva14,BasiPana16,EstrFrid16,RodrSanc17,Tren19,MercUrib21}, and under additional regularity assumptions those in~\cite{YangDing18}.

As the main contribution we develop a feedback controller, which is designed to achieve satisfaction of a particular control objective.
While the recently proposed funnel controller~\cite{BergIlch21} achieves asymptotic exact tracking, the controller in the present article achieves \emph{exact tracking in predefined finite time}.
The latter means that the output~$y$ of a system is forced to approach a given reference~$y_{\rm ref}$, and $\lim_{t \to T} y(t) = y_{\rm ref}(T)$ for a predefined final time~$T$, cf. Figure~\ref{Fig:ControlObjective}.
The result in the present article closes a gap in the existing theory of \emph{funnel control}, and is formulated in Theorem~\ref{Thm:Exact-tracking-in-finte-time}.
We rigorously prove feasibility of the proposed controller.
To this end, we extend and generalize the existing feasibility proof~\cite{BergIlch21}.
Specifically, the proof in~\cite{BergIlch21} extensively uses a growth condition on the involved ``funnel function''~$\vp$, while in the present article the respective function is unbounded, it even has a pole.
Since the proof is quite technical, it is relegated to the appendix.
\blfootnote{
\textbf{Nomenclature.}
$[a,b], [a,b), (a,b)$ is a closed, half-open, and open interval for~$a,b \in \R$,~$a < b$;
$\rp := [0, \infty)$;
$ \langle \cdot, \cdot \rangle $ is the inner product in~$\R^n$;
$\| x \| := \sqrt{ \langle x, x \rangle }$ for $x \in \R^n$;
$\Gl_n(\R)$ is the group of invertible $\R^{n\times n}$ matrices;
for $I \subseteq \R$ an interval
$\cL_{\rm loc}^\infty (I ; \R^p)$ is the set of locally essentially bounded functions $ f: I \to \R^p$;
$\cL^\infty (I ; \R^p)$ is the set of essentially bounded functions $ f: I \to \R^p$ ;
$\| f \|_{\infty} := \esssup_{t \in I} \|f(t)\| $ norm of $f \in \cL^\infty(I ; \R^p)$;
$\cW^{k,\infty}(I ; \R^p)$ is the set of $k$-times weakly differentiable functions $ f : I \to \R^p$ such that $f,\dot f,\ldots,f^{(k)} \in \cL^\infty(I ; \R^p)$;
$\cC^k( I ; \R^p) $ is the set of $k$-times continuously differentiable functions $f : I \to \R^p$, $\cC(I; \R^p) = \cC^0(I; \R^p)$;
$f|_{J}$ is the restriction of $f : I \to \R^n$ to $J \subseteq I$, $I,J$ intervals. 
}
%
\section{Control objective, system class, feedback law}
We state the problem under consideration, introduce the class of systems to be controlled, and define the feedback law.
To convey the basic idea, we briefly give the general framework and then present the individual components in detail in Sections~\ref{Sec:SystemClass} and~\ref{Sec:FeedbackLaw}.
We consider multi-input multi-output $r^{\rm th}$-order functional differential equations
\begin{equation} \label{eq:system}
\begin{aligned}
y^{(r)}(t) &= f\big( d(t), \textbf{T}(y,\dot y,\ldots, y^{(r-1)})(t), u(t) \big), \\
y|_{[-\sigma,0]} &= y^0 \in \cC^{r-1}([-\sigma,0] ; \R^m),
\end{aligned}
\end{equation}
with bounded \emph{unknown} disturbance~$d$, \emph{unknown} nonlinear function~$f$ and \emph{unknown} operator~$\oT$, the latter are characterized in Definitions~\ref{Def:high-gain-prop} and~\ref{Def:OP-class} below.
If $\sigma > 0$, then the initial value is given via the initial trajectory~$y^0$; if $\sigma=0$, then the initial value is given by $(y(0),\dot y(0),\ldots,y^{(r-1)}(0)) \in \R^{rm}$.
Beside typical physical phenomena such as, e.g., hysteresis effects, the operator~$\oT$ can also model delay elements, cf.~\cite[Sec.~4.4]{IlchRyan02b}. If delays are involved, $\sigma > 0$ corresponds to the initial delay.
Note that the input~$u$ and the output~$y$ have the same dimension~$m \in \N$. \\
For a control function $u \in \cL^\infty_{\rm loc}(\rp, \R^m)$, system~\eqref{eq:system} has a solution in 
the sense of \textit{Carath\'{e}odory}, meaning a function
$x : [-\sigma,\omega) \to \R^{rm}$, $\omega > 0$, 
with $x|_{[-\sigma,0]} = (y^0,\dot y^0,\ldots,(y^0)^{(r-1)})$
such that $x\vert_{[0,\omega)}$ is absolutely continuous and 
satisfies $\dot x_i(t) = x_{i+1}(t)$ for $i=1,\ldots,r-1$, and $\dot x_r(t) = f(d(t),\oT(x(t)), u(t))$ (which corresponds to~\eqref{eq:system}) for almost all~$t\in[0,\omega)$.
A solution $x$ is said to be \textit{maximal}, if it does not have a right extension which is also a solution.
\subsection{Control objective} 
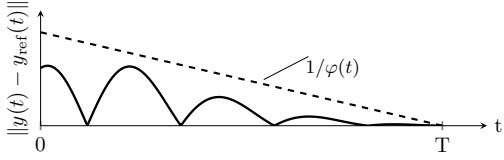
\begin{figure}
\centering
\begin{tikzpicture}[scale=0.8]
\begin{axis}[
axis lines = left, ylabel near ticks, ylabel = \( \|y(t)-y_{\rm ref}(t)\| \),
xmin=0, xmax=10,
ymin=0, ymax=10,
scale=1,
xtick={0,9},
x tick style = {thick},
xticklabels={0,T},
xlabel = t,
every axis x label/.style={
    at={(ticklabel* cs:1.0)},
    anchor=west,
},
ytick=\empty,
yticklabels={},
height = 3.5cm,
width = 9cm,
]
\addplot[dashed, domain=0:9,samples=10, line width=1.0pt] {0.9*(9-x)} ;
\addplot[domain=0:9,samples=300 , line width=1.1pt] {abs(5*exp(-0.5*x)*(x+1)*cos(deg(1.5*x))*cos(deg(pi/18*x)))} ;
    \draw (axis cs:5,4) -- node[above=0.4em,right=0.2cm]{\footnotesize{$ 1/\vp(t)$}} (axis cs:6,6);
\end{axis}
\end{tikzpicture}
\caption{Schematic illustration of the control objective~\eqref{eq:control-objective}.}
\label{Fig:ControlObjective}
\end{figure}

We aim to design a feedback controller, which achieves exact reference tracking in the following sense.
For a given reference trajectory $y_{\rm ref} \in \cW^{r,\infty}([0,T) ; \R^m)$ and a predefined final time~$T > 0$, the output~$y$ of system~\eqref{eq:system} approaches the reference within the interval~$[0,T)$, and coincides with the reference as $t \to T$, i.e., for $e(\cdot) := y(\cdot) - y_{\rm ref}(\cdot)$
\begin{subequations} \label{eq:control-objective}
\begin{equation} 
\begin{aligned}
 \forall \, i=0,\ldots,r-1 \,  : \
 \lim_{t \to T} \|e^{(i)}(t)\| = 0,
\end{aligned}
\end{equation}
where~$e^{(i)}(\cdot)$ denotes the~$i^{\rm th}$ derivative of~$e(\cdot)$.
Moreover, in the transient phase for $t \in [0,T)$ the error evolves within the so called ``performance funnel", i.e.,
\begin{equation}
\begin{aligned}
 \! \! \! (t,e(t)) \!  \in \!  \cF_\vp \!  := \! \setdef{ \! (t,e) \! \in \![0,T) \!  \times \! \R^m \!}{ \vp(t) \|e\| < 1 },
\end{aligned}
\end{equation}
\end{subequations}
where $\vp$ is a boundary function defined in~\eqref{def:vp}. 
The control objective is illustrated in Figure~\ref{Fig:ControlObjective}. 

\subsection{System class} \label{Sec:SystemClass}
To introduce the system class under consideration, we first provide some necessary definitions.
To characterize the class of admissible nonlinearities~$f$ in system~\eqref{eq:system}, we recall the definition of the ``high-gain property" from~\cite[Sec.~1.2]{BergIlch21}.
\begin{Definition}  \label{Def:high-gain-prop}
For~$p,q,m \in \N$ a function $f \in \cC(\R^p \times \R^q \times \R^m ; \R^m)$ satisfies the \textit{high-gain property}, if there exists~$\rho \in (0,1)$ such that, for every compact~$K_p \subset \R^p$ and compact~$K_q \subset \R^q$ the continuous function
\begin{equation*}
\begin{aligned} 
\chi : \R & \to \R, \\
s &\mapsto \min \setdef{ \langle v, f(\delta,z,-s v \rangle }{ \begin{array}{l}
(\delta,z) \in K_p \times K_q,  \\
v \in \R^m, \rho \le \| v \| \le 1 
\end{array} }
\end{aligned}
\end{equation*}
is such that $\sup_{s \in \R} \chi(s) = \infty$.
\end{Definition}
In Remark~\ref{Rem:Comments-on-parameters} we discuss the high-gain property in detail.
The operator $\oT$ in~\eqref{eq:system} belongs to the operator class defined below.
This definition is taken from~\cite[Sec.~1.2]{BergIlch21}.
\begin{Definition} \label{Def:OP-class} 
If for $n,q \in  \N$ and $\sigma \ge 0$ the operator $\textbf{T} : \cC([-\sigma,\infty) ; \R^{n}) \to \cL_{\rm loc}^\infty(\rp ; \R^q)$ satisfies
\begin{enumerate}[label = (\alph{enumi}), ref=(\alph{enumi}),leftmargin=*]
\setlength \itemsep{0cm}
\item $\textbf{T}$ maps bounded trajectories to bounded trajectories, i.e., for all $c_1 > 0$, there exists $c_2>0$ such that for all 
$\xi \in \cC([-\sigma,\infty) ; \R^n)$,
\[
\sup_{t \in [-\sigma,\infty)} \| \xi(t) \| \le c_1 \ \Rightarrow \ \sup_{t \in [0,\infty)} \| \textbf{T}(\xi)(t) \| \le c_2,
\]
\item $\textbf{T}$ is causal, this means, for all $t \ge 0$ and all functions $\zeta, \xi \in \cC([-\sigma,\infty) ; \R^n)$, 
\[
\zeta |_{[-\sigma,t)} = \xi|_{[-\sigma,t)} \ \Rightarrow \ \textbf{T}(\zeta)|_{[0,t)} \overset{a.a.}{=} \textbf{T}(\xi)|_{[0,t)},
\]
\item $\textbf{T}$ is locally Lipschitz continuous in the following sense: for all $t \ge 0 $ and all $\xi \in \cC([-\sigma,t] ; \R^n)$ 
there exist $\Delta, \delta, c > 0$ such that for all 
$\zeta_1, \zeta_2 \in \cC([-\sigma,\infty) ; \R^n)$ with $\zeta_1|_{[-\sigma,t]} = \xi $, $\zeta_2|_{[-\sigma,t]} = \xi $ 
and $\|\zeta_1(s) - \xi(t)\| < \delta$,  $\|\zeta_2(s) - \xi(t)\| < \delta$ for all $s \in [t,t+\Delta]$ we have
\begin{align*}
\esssup_{s \in [t,t+\Delta]} & \| \textbf{T}(\zeta_1)(s) - \textbf{T}(\zeta_2)(s) \| \\ 
& \le c \ \rm{sup}_{s \in [t,t+\Delta]} \| \zeta_1(s) - \zeta_2(s)\| ,
\end{align*}
\end{enumerate}
then the operator~$\textbf{T}$ belongs to the operator class~$\cT_\sigma^{n,q}$.
\end{Definition}
With the definitions so far, we may introduce the system class under consideration, which is the same as in~\cite{BergIlch21}.
\begin{Definition}
For~$m,r \in \N$ a system~\eqref{eq:system} is said to belong to the system class~$\cN^{m,r}$, if for some~$p,q \in \N$ the ``disturbance" is bounded, i.e., $d \in \cL^\infty(\rp ; \R^p)$, the function $f \in \cC(\R^p \times \R^q \times \R^m ; \R^m)$ satisfies the high-gain property from Definition~\ref{Def:high-gain-prop} and for~$\sigma \ge 0$ the operator~$\textbf{T}$ belongs to~$\cT^{rm,q}_{\sigma}$; we write~$(d,f,\textbf{T}) \in \cN^{m,r}$.
\end{Definition}
\begin{Remark}
For $n \in \N$,  consider a state-space model
\begin{equation} \label{eq:state-space}
\begin{aligned}
\dot x(t) &= \tilde f(x(t)) + \tilde g(x(t)) u(t),  \quad x(0) = x^0 \in \R^n, \\
y(t) &= \tilde h(x(t)),
\end{aligned}
\end{equation}
where $y(t) \in \R^m$ is the output, and for $m \le n$,
$\tilde f: \R^n \to \R^n$, $\tilde g: \R^n \to \R^{n \times m}$, $\tilde h : \R^n \to \R^m$ sufficiently smooth.
Assume this system has relative degree~$(r_1,\ldots,r_m) = (r,\ldots,r) \in \N^m$ for~$r \in \N$ at~$x_0 \in \R^n$, i.e., there exists a neighbourhood $U \subseteq \R^n$ of~$x_0$ such that
\begin{align*}
 \forall \, k = 0,\ldots,r-2\, \forall z \in U: \  & (L_{\tilde g}^{} L_{\tilde f}^{k} \tilde h)(z) = 0, \\ 
\text{and} \ \forall z \in U : \    \gamma(z) :=  &(L_{\tilde g} L_{\tilde f}^{r-1} \tilde h)(z)  \in \Gl_m(\R),
\end{align*}
where $(L_{\tilde f} \tilde h)(z) :=\tilde  h'(z) \cdot \tilde  f(z)$ denotes the \emph{Lie derivative} of~$\tilde h$ along~$\tilde f$.
Then, by~\cite[Prop.~5.1.2]{Isid95}, there exists a coordinate transformation~$\Phi: U \to W$, $W \subseteq \R^n$ open, transforming~\eqref{eq:state-space} into Byrnes-Isidori form
\begin{align*}
\dot \xi_i(t) &= \xi_{i+1}(t) \in \R^m, \quad i=1,\ldots,r-1, \\
\dot \xi_r(t) &= (L_{\tilde f}^r \tilde h)(\Phi^{-1}(\xi(t),\eta(t))) + \gamma(\Phi^{-1}(\xi(t),\eta(t))) u(t), \\
\dot \eta(t) & = q(\xi(t),\eta(t)) ,
\end{align*}
where $\xi_1(t) = y(t) \in \R^m$ is the original output, and $\eta$ denotes the internal dynamics.
If~$\gamma(\cdot)$ is sign definite, then, for appropriate~$f,\oT$, the state-space representation~\eqref{eq:state-space} is  locally equivalent to~\eqref{eq:system}.
In~\cite{ByrnIsid91a,schwartz1999global} sufficient conditions for a global transformation are formulated in terms of differential geometric properties.
In the present article, structural conditions are formulated in terms of the high-gain property (Definition~\ref{Def:high-gain-prop}) and the operator class (Definition~\ref{Def:OP-class}).
Systems~\eqref{eq:system} are restricted to $r_i = r$ for all $i=1,\ldots,m$.
It is possible to generalize the proposed controller to systems
\begin{equation*}
        (y_1^{(r_1)}(t),\ldots,y_m^{(r_m)}(t)) ^\top
        = f(t,\Lambda(y)(t), \oT(\Lambda(y))(t), u(t))
\end{equation*}
where $\Lambda(y)(t) := (y_1(t),\dot y_1(t),\ldots,y_m^{(r_m-1)}(t))$.
However, since such systems massively increase complexity of notation, we 
restrict ourselves to systems~\eqref{eq:system}.
\end{Remark}

\subsection{Feedback law} \label{Sec:FeedbackLaw}
We formulate the feedback law, which achieves the control objective~\eqref{eq:control-objective}.
The two main ingredients are the prescribed final time~$T > 0$, and the error boundary, i.e., the funnel function~$\vp$.
To establish the controller, we introduce the following control parameters. 
Choose the final time $T > 0$, some~$c>0$,
and the funnel function
\begin{subequations} \label{def:control_functions}
\begin{equation} \label{def:vp}
 \vp(t) = \frac{1}{c} \frac{1}{T-t}, \quad t \in [0,T).
\end{equation}
Note that~$\lim_{t \to T} \vp(t) = \infty$. 
For~$c>0$ in~\eqref{def:vp} choose
\begin{equation} \label{def:alp}
\alpha \in \cC^{r-1}( [0,1) ; [c(r+1) , \infty)) \ \text{a bijection},
\end{equation}
and further choose
\begin{equation} \label{def:surjection}
N \in \cC(\rp ; \R) \ \text{a surjection}.
\end{equation}
\end{subequations}
In Remark~\ref{Rem:Comments-on-parameters} we comment on the control parameters defined in~\eqref{def:control_functions}.
%
We set~$e^{(i)}(\cdot) := y^{(i)}(\cdot) - y_{\rm ref}^{(i)}(\cdot)$ for $i=0,\ldots,r-1$, and recursively define for $\gamma_0 = 0$, and~$k=1,\ldots,r$ with~$\alpha(\cdot)$ from~\eqref{def:alp} 
the functions
\begin{subequations} \label{def:control_variables}
\begin{align}
e_k(t) &= \vp(t) e^{(k-1)}(t) + \vp(t) \sum_{i=1}^{k-1} \gamma_{i}^{(k-1-i)}(t), \label{def:e_k} \\
\gamma_k(t) &=  \alpha(\|e_k(t)\|^2) e_k(t).  \label{def:gam_k}
\end{align}
\end{subequations}
Then, with the functions introduced in~\eqref{def:control_functions},~\eqref{def:control_variables} we define the feedback law $u : \rp \to \R^m$ by
\begin{equation} \label{def:u}
u(t) := (N \circ \alpha)(\|e_r(t)\|^2) \, e_r(t) .
\end{equation}
\begin{Remark} \label{Rem:Comments-on-parameters}
We comment on the control parameters defined in~\eqref{def:control_functions}, and on the high-gain property.
\begin{enumerate}[label = \roman{enumi}), leftmargin=1.4em]
\item The high-gain property in Definition~\ref{Def:high-gain-prop} is essential to achieve the control objective~\eqref{eq:control-objective}.
If a large input is applied, the system has to react appropriately, i.e., if the error is close to the funnel boundary, a large input results in a ``fast'' response of the system.
For a more detailed discussion and equivalent conditions of the high-gain property we refer to~\cite[Rem.~1.3~\&~1.4]{BergIlch21}.
\item Compared to the funnel controller proposed in, e.g., \cite{IlchRyan02b,BergIlch21}, the explicit choice of~$\vp(\cdot)$ in~\eqref{def:vp} seems restrictive.
Anticipating the initial conditions~\eqref{eq:initial} in Theorem~\ref{Thm:Exact-tracking-in-finte-time}, this choice of~$\vp(\cdot)$ reflects the intuition that the shorter the final time~$T$ is chosen, the better the initial guess has to be. 
\item The bijection~$\alpha(\cdot)$ is responsible for the \emph{high-gain}, i.e., the smaller the distance between the error and the funnel boundary is, the larger the input values are. A typical choice is $\alpha(s) = c(r+1)/(1-s)$.
\item The parameter~$c>0$ in~\eqref{def:vp} links the funnel function~$\vp(\cdot)$ to the gain function~$\alpha(\cdot)$ in~\eqref{def:alp}.
The larger the value~$c>0$, the larger the lower bound of~$\alpha(\cdot)$, i.e., small tracking errors result in higher input values.
From the perspective of the initial conditions~\eqref{eq:initial}, the parameter~$c>0$ can be used to satisfy these. Given a final time $T>0$ and initial values, the inequalities in~\eqref{eq:initial} can be utilized to find an appropriate $c>0$.
So in view of the initial conditions, the parameter $c>0$ and the final time $T$ have an intuitive relation: the shorter $T$ is, the larger $c$ must be.
Moreover, the larger the initial error is, the large~$c$ must be.
\item The surjection~$N(\cdot)$ from~\eqref{def:surjection} accounts for possible unknown control directions. A feasible choice is, e.g, $N(s) = s \sin(s)$.
If the control direction is known, e.g., $y^{(r)}(t) = f(d(t),\oT(y,\ldots,y^{(r-1)}(t)) + \Gamma u(t)$ with $\Gamma > 0$ (or $\Gamma<0$), then the simple choice $N(s) = -s$ (or $N(s) = s$) is feasible.
For detailed comments on the surjection~$N(\cdot)$ and possible further simplifications we refer to~\cite[Rem.~1.8]{BergIlch21}.
\end{enumerate}
\end{Remark}
\begin{Remark} \label{Rem:Gam_k}
We comment the computation of~\eqref{def:control_variables}. 
Define the set
$
\cD_{0} := \setdef{ \zeta \in \R^m }{ \| \zeta \| < 1 } , 
$
and, for~$\alpha(\cdot)$ from~\eqref{def:alp}, the function~$\Gamma_{0}: \cD_{0} \to \R^m$ by
$
\Gamma_{0}(t,\zeta) := \alpha(\|\zeta\|^2) \, \zeta.
$
For~$k=1,\ldots,r-1$ recursively define the sets~$\cD_{k}$ and the functions $\Gamma_{k} : [0,T) \times \cD_{k} \to \R^m$ by
\begin{align} \label{def:Gam_k}
\cD_k &:= \underbrace{\cD_0 \times \cdots \times \cD_0}_{k-\rm times} \times \R^m, \nonumber \\
& \Gamma_{k}(t,\zeta_1,\ldots,\zeta_{k+1}) := \frac{\partial \Gamma_{k-1}(t,\zeta_1,\ldots,\zeta_{k})}{\partial t} \\
& \!+\! \sum_{i=1}^{k} \frac{\partial \Gamma_{k-1}(t,\zeta_1,\ldots,\zeta_{k})}{\partial \zeta_i} 
\Big( \! \vp(t) \left( c \zeta_i \! -\! \Gamma_{0}(\zeta_i) \right) \! +\! \zeta_{i+1} \!\Big). \nonumber
\end{align}
Then, with~$e_k(\cdot)$ from~\eqref{def:e_k} and $\gamma_k(\cdot)$ from~\eqref{def:gam_k} we obtain 
$
\gamma_k^{(j)}(t) = \Gamma_j \big( t,e_k(t),\ldots,e_{k+j}(t) \big), 
$
for $0 \le j \le r-k$,
which can be seen via a induction over~$k$ using~\eqref{def:control_variables}.
\end{Remark}
%
Due to the recursion~\eqref{def:control_variables}, the controller~\eqref{def:u} is not as simple to implement as the controller in \cite[Eq.~(9)]{BergIlch21}.
However, given~\eqref{def:Gam_k}, the calculation of the required expressions can be done completely algorithmically.
%

\section{Main result} \label{Sec:Main-result}
This section contains the main result. 
To phrase it, the application of the controller \eqref{def:u} to a system \eqref{eq:system} yields a closed-loop initial value problem that has a solution; the input and output signals are bounded and in particular, the controller achieves exact output tracking in predefined finite time with prescribed behaviour of the tracking error.
\begin{Theorem} \label{Thm:Exact-tracking-in-finte-time}
For $m,r \in \N$ consider a system~\eqref{eq:system} with $(d,f,\textbf{T}) \in \cN^{m,r}$ and initial data $y^0 \in \cC^{r-1}([-\sigma,0] ; \R^m)$.
Let~$T>0$ and $y_{\rm ref} \in \cW^{r,\infty}([0,T) ; \R^m)$.
Assume that with the control parameters in~\eqref{def:control_functions} the auxiliary control variables $e_k(\cdot)$ in~\eqref{def:control_variables}
satisfy the initial conditions
\begin{equation} \label{eq:initial}
\forall\, k=1,\ldots,r \, : \ \|e_k(0)\| < 1.
\end{equation}
Then, the funnel controller~\eqref{def:u} applied to~\eqref{eq:system} yields an initial value problem, which has 
a solution and every maximal solution~$y: [-\sigma, \omega) \to \R^m$ satisfies
\begin{enumerate}[label=\roman*)]
\item $\omega = T$,
\item $u \in \cL^\infty([0,T) ; \R^m)$, $y \in \cW^{r,\infty}([-\sigma, T) ; \R^m)$,
\item the tracking error~$e(t) = y(t) - y_{\rm ref}(t)$ evolves within the performance funnel~$\cF_\vp$, i.e., 
\begin{equation*}
 \forall \, t \in [0,T) \,: \ \vp(t) \| e(t) \| < 1,
\end{equation*}
\item the tracking of the reference and its derivatives is exact at $t=T$, i.e.,
\begin{equation*}
\forall\, i=0,\ldots,r-1 \, : \ \lim_{t \to T} \|e^{(i)}(t)\| = 0.
\end{equation*}
\end{enumerate}
\end{Theorem}
The proof is relegated to the appendix.
Note that, since the system class~$\cN^{m,r}$ encompasses the systems under consideration in \cite{Andr08,Leva14,BasiPana16,EstrFrid16,RodrSanc17,Tren19,MercUrib21}, and under additional regularity assumptions those in~\cite{YangDing18}, the proposed feedback law~\eqref{def:u}, assuming availability of the first $r-1$ output derivatives, achieves the control objectives formulated in those references with prescribed behaviour of the error and within predefined finite time.
\begin{Remark}
Since at the first glance the control law~\eqref{def:u} is very similar to the controller proposed in~\cite{BergIlch21}, we emphasize some differences.
\begin{enumerate}[label = \roman{enumi}), ref=\roman{enumi}), leftmargin=1.4em]
\item As highlighted in, e.g., ~\cite{IlchRyan09,BergIlch21}, some care is required when showing boundedness of the involved signals, since the bijection~$\alpha(\cdot)$ may introduce a singularity.
Moreover, in the present context, expressions involving the unbounded funnel function~$\vp(\cdot)$ demand particularly high attention, 
cf. \emph{Steps two} and \emph{three} in the proof.
\item
A careful inspection of the proof of~\cite[Thm.~1.9]{BergIlch21} yields, that the following growth condition 
\begin{equation*}
\exists \, d > 0  : \ |\dot \phi(t) | \le d \big( 1+\phi(t) \big) \text{ for almost all } t \ge 0
\end{equation*}
on the funnel function~$\phi$ is crucial.
It prevents a ``blow up" in finite time, i.e., $\phi(\cdot)$ is bounded on any compact interval.
With this, however, exact tracking in finite time via funnel control is impossible.
Contrary, the funnel functions~$\vp(\cdot)$ in~\eqref{def:vp} do not satisfy this growth condition.
Hence, the respective steps in the proof of
 \cite[Thm.~1.9]{BergIlch21} are not valid in the present analysis.
\item 
In order to show boundedness of the involved error signals, novel techniques have been developed in the proof of Theorem~\ref{Thm:Exact-tracking-in-finte-time}. 
In particular \emph{Steps two, three} and \emph{four} contain innovations not found in the existing works on high-gain feedback control.
\item The conclusion drawn in \emph{Step eight}, namely that the tracking error is zero at $t=T$, is only possible with the results derived in \emph{Step three}.
\end{enumerate}
\end{Remark}
\begin{Remark}
Assertion~$i)$ in Theorem~\ref{Thm:Exact-tracking-in-finte-time}, namely~$[0,T)$ being the maximal solution interval, naturally raises the question of a global solution in time.
If the system's equations~\eqref{eq:system} are available and~$y_{\rm ref}(\cdot)$ is defined on~$\rp$, then $y(t) - y_{\rm ref}(t) \equiv 0$ for all time~$t \ge T$ can be achieved by asking the reference to satisfy~\eqref{eq:system} for~$ t \ge T$, with~$u(\cdot) \equiv 0$, and ``initial" conditions~$y_{\rm ref}(T) = y(T), \dot y_{\rm ref}(T) = \dot y(T),\ldots, y_{\rm ref}^{(r-1)}(T) = y^{(r-1)}(T) $.
\end{Remark}

\begin{Remark} \label{Rem:practical-exact}
For any given $\ve > 0$ there exists a time $T_\ve < T$ such that each of the first $r-1$ derivatives of the error $e(\cdot) = y(\cdot) - y_{\rm ref}(\cdot)$ can be bounded by $\ve$ for all $t \in [T_\ve , T)$, this is, for all $k=0,\ldots,r-1$
and all $t \in [T_\ve, T)$ we have $\|e^{(k)}(t)\| < \ve$.
This property is relevant, e.g., if during a docking manoeuvre the demanded accuracy changes.
\end{Remark}

\section{Numerical examples} \label{Sec:Simulation}
We present two numerical examples.
In the first simulation we consider a docking maneuver as an application.
The second simulation illustrates how the choice of~$\vp$ affects the maximal control input.
\subsection{Docking maneuver}
\begin{figure}[h]
\begin{center}
\adjustbox{width=0.9\linewidth}{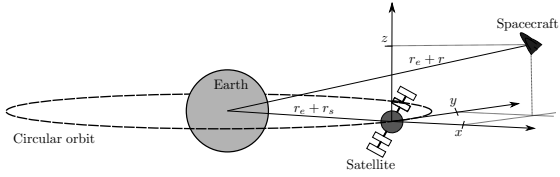}
\caption{Hill's local-vertical-horizontal coordinate frame.}
\label{Fig:Hills_Frame}
\end{center}
\end{figure}

\begin{figure}
\centering
\includegraphics[width=0.76\linewidth]{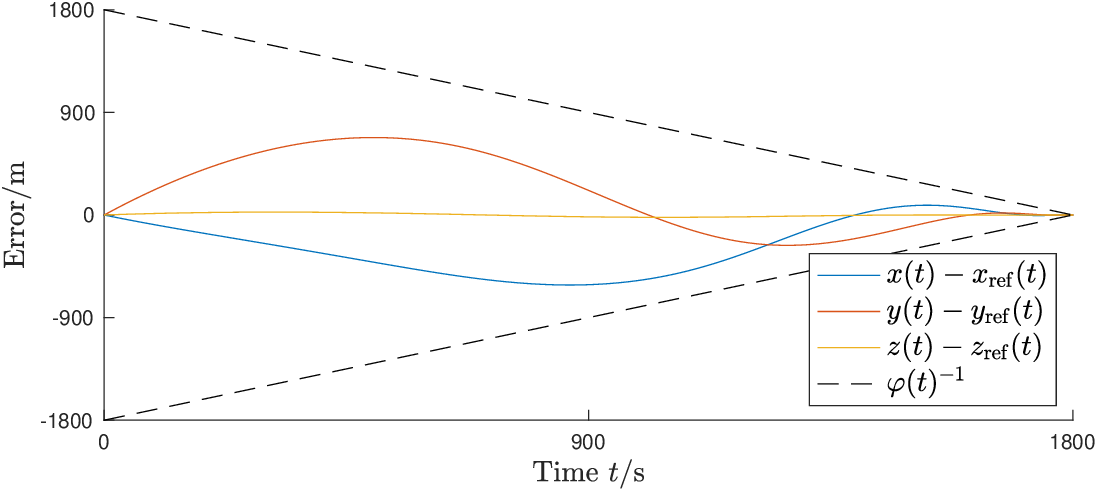}
\caption{Errors and funnel.}
\label{Fig:Errors_and_funnel}
\end{figure}

\begin{figure}
\centering
\includegraphics[width=0.73\linewidth]{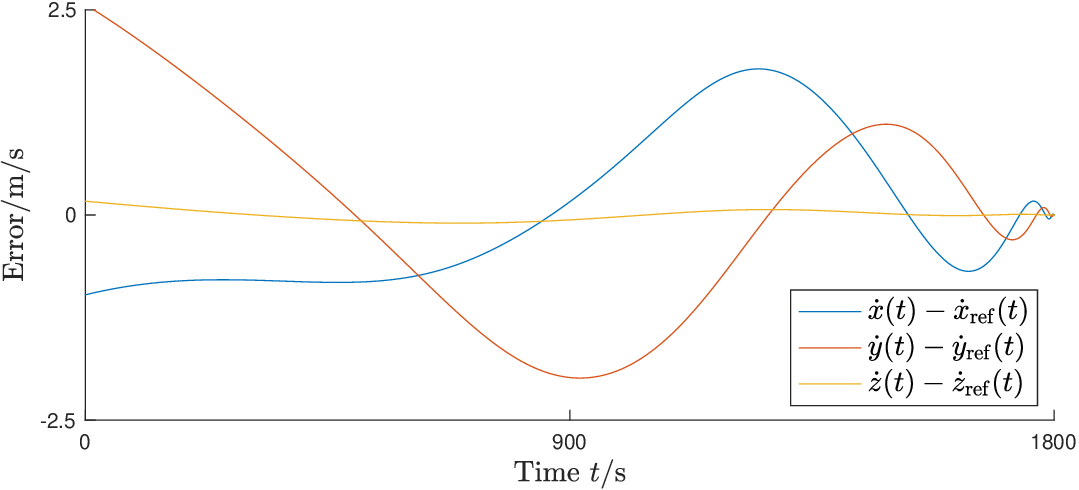}
\caption{Error of the respective velocities.}
\label{Fig:Errors_velocities}
\end{figure}

\begin{figure}
     \centering
     \begin{subfigure}[b]{\linewidth}
         \centering
         \includegraphics[width=0.73\linewidth]{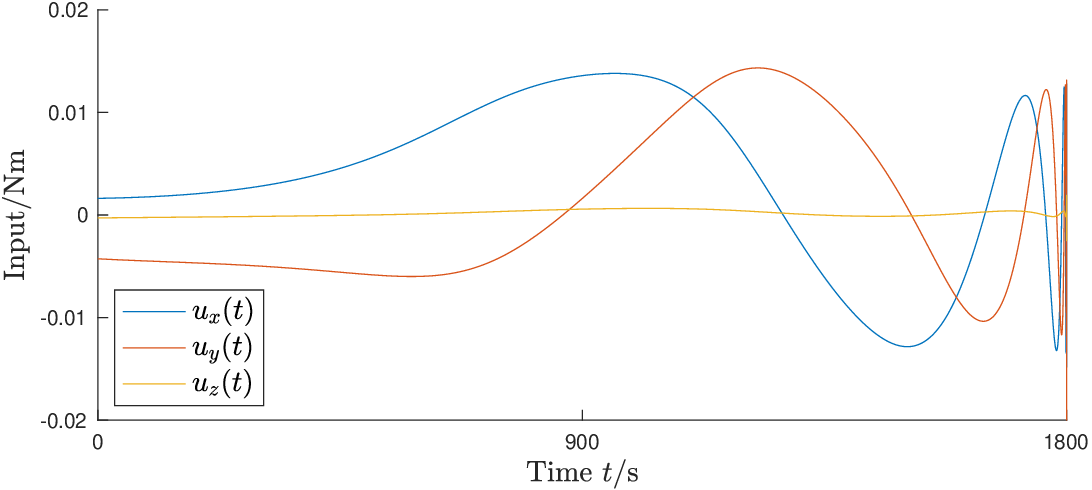}
			\caption{Control input.}
			\label{Fig:Control_detail}
     \end{subfigure}
     \hfill
          \begin{subfigure}[b]{\linewidth}
         \centering
         \includegraphics[width=0.73\linewidth]{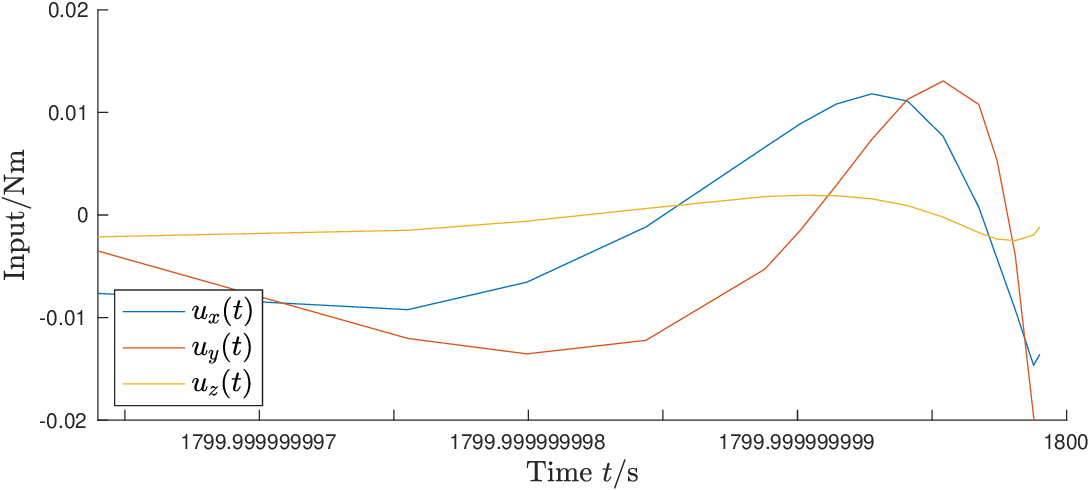}
			\caption{Control input in the very last moments.}
			\label{Fig:Control_detail_end}
     \end{subfigure}
     \caption{Control input.}
     \label{Fig:Control}
\end{figure}

As an exemplary application we simulate docking of two spaceships.
Consider a passive space station in a circular orbit, and a chasing active spacecraft. 
We assume the passive space station to be on a constant altitude~$r_s$ with constant angular velocity~$\omega = \sqrt{\mu / (r_e + r_s)^3}$, where~$\mu \approx 3.986 \cdot 10^{14} \, \rm{m^3/s^2}$ is the standard gravitational parameter,
 and~$r_e = 6378137 \, \rm{m}$ the radius of the earth.
To analyze the motion of the spacecraft, we use Hill's local-vertical-local-horizontal coordinate frame~\cite{Hill78}, see Figure~\ref{Fig:Hills_Frame}.
Within this frame we use the commonly used Clohessy-Wiltshire model for satellite rendezvous~\cite{ClohWilt60}, also elaborated on in~\cite{Klue99}.
Let~$r(t)$ denote the altitude of the chasing spacecraft at time~$t$. 
Then, for~$x(\cdot) = r(\cdot) - r_s(\cdot)$ the component of relative distance along the radial direction, 
$y(\cdot)$ the downtrack component along satellite's circular orbit, 
and~$z(\cdot)$ the distance component along the satellite's angular momentum, 
and setting~$ \zeta(\cdot) := (\zeta_1(\cdot),\zeta_2(\cdot),\zeta_3(\cdot))^\top := (x(\cdot),y(\cdot),z(\cdot))^\top$ we obtain the Clohessy-Wiltshire equations
\begin{equation*}
\begin{scriptsize}
\begin{aligned}
\ddot \zeta_1(t) &= 3 \omega^2 \zeta_1(t) + 2\omega \dot \zeta_2(t) + u_x(t) +d_x(t), \\
\ddot \zeta_2(t) &= -2\omega \dot \zeta_1(t) + u_y(t) + d_y(t), \\
\ddot \zeta_3(t) &= -\omega^2 \zeta_3(t) + u_z(t) + d_z(t).
\end{aligned}
\end{scriptsize}
\end{equation*}
Setting 
$f: \R_{\ge 0} \R^3 \times \times \R^3 \to \R^3$, $(t,\xi, \nu) \mapsto (d_x(t) + 3 \omega^2 \xi_1 + 2\omega \nu_2, d_y(t)  -2 \omega \nu_1, d_z(t)-\omega^2 \xi_3)^\top $,
with $B = I_3 \in \R^{3 \times 3}$ the equations of motion above with output~$ \zeta(\cdot)$ can compactly be written as
$\ddot \zeta(t) = f(t,\zeta(t), \dot \zeta(t)) + Bu(t)$, 
which is a system belonging to~$\cN^{3,2}$ with disturbance.
For simulation purposes we choose $r_s =415000 \, \rm m $ ($\approx$~altitude ISS), which yields $\omega \approx 0.00113 \, \rm{s^{-1}}$ corresponding to an orbital period of approximately $93$ minutes.
Since we simulate a docking maneuver, we choose the reference $\zeta_{\rm ref}(\cdot)$ such that $\zeta(T) = (0,0,0)^\top$ in Hill's coordinate frame.
Let
$
\zeta_{\rm ref}(t) = \zeta(0) ( 1- \sin (\frac{\pi}{2} \tfrac{t}{T} ) ),
$
guiding the spacecraft smoothly to the satellite, with $\zeta_{\rm ref}(T) = \dot \zeta_{\rm ref}(T) = (0,0,0)^\top$.
We take the initial conditions from~\cite{Klue99} $x(0)=-y(0) = 1000\, \rm m$, and additionally~$z(0)=250 \, \rm m$;
and~$\dot x(0) = - 0.1 \, \rm{ms^{-1}}$, $\dot y(0)= 1.69 \, \rm{ms^{-1}}$, and~$\dot z(0)=-0.05 \, \rm{ms^{-1}}$.
As docking time we choose~$T = 1800 \, \rm s$.
As control parameters we choose $N : s \mapsto -s $, and $\alpha: s \mapsto (r+1)c/(1-s)$. 
The disturbance is $d_i(t) = 0.05\sin(t)\cos(0.7t)$, $i=x,y,z$.
With $c = 1$ the initial conditions~\eqref{eq:initial} are satisfied.
Since $\lim_{t \to T} \vp(t) = \infty$, simulation is only possible for $[0,t_{\rm max}]$ with $t_{\rm max} < T$.
Since $\vp(t) \|e(t)\| < 1$ for all $t \in [0,T)$, in virtue of Remark~\ref{Rem:practical-exact} a value $\textsf{eps}$ can be chosen such that a certain upper bound of the spatial error at final time $t_{\rm max}$ is guaranteed, i.e.,
$
\| e(t_{\rm max}) \| < \frac{1}{\vp(t_{\rm max})} = c(T-t_{\rm max}) \le \textsf{eps},
$
from which we obtain $t_{\rm max} \ge T - \textsf{eps}/c$.
Here we choose~$\textsf{eps} = 10^{-10} \, \rm m$, which means a spatial accuracy of \AA ngstr\"om (range of size of atoms).
This seems to be a unnecessary high accuracy since in real applications the required rendezvous distance is about centimetres, then magnetic docking structures become active; however, if these fail unexpectedly, the feedback control still is capable to perform a docking maneuver.
Simulations have been performed in \textsc{Matlab} (solver: \textsf{ode15s}, \textsf{AbsTol=RelTol=$10^{-12}$}).
Figure~\ref{Fig:Errors_and_funnel} shows  that the docking maneuver is successful within the predefined finite time~$T$, and the errors evolve within the prescribed boundary.
Figure~\ref{Fig:Errors_velocities} shows the errors of the velocities, where the disturbance can be seen as fast oscillations. 
As expected from Theorem~\ref{Thm:Exact-tracking-in-finte-time}, the errors of the velocities tend to zero for~$t \to T$.
The control input is depicted in Figure~\ref{Fig:Control}.
Figure~\ref{Fig:Control_detail_end} shows the control input in the very last moments before docking, where the largest input signals are generated, which may not be needed if the docking tools are activated and work as intended.
\subsection{Control effort and funnel function} \label{Sec:u_c}
To illustrate how the choice of~$\vp$ in~\eqref{def:vp} influences the control input, we consider a chain of integrators
\begin{equation*}
\ddot y(t) = u(t), \quad y(0)  = 1, \ \dot y(0)  = 0,
\end{equation*}
and perform stabilization, i.e., $y_{\rm ref}(t) \equiv 0$, with exact value at final time~$T = 1$. 
To influence the shape of~$\vp$, the parameter~$c$ is varied with $c_k= 2 + 0.1 k$ for $k=0,\ldots,200$, while keeping the final time~$T$ constant.
\begin{figure}
\centering
\includegraphics[width=0.73\linewidth]{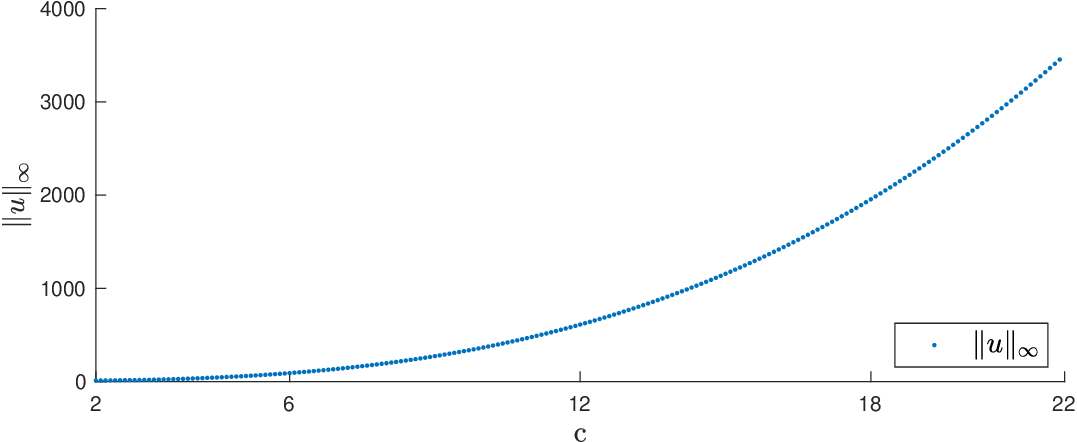}
\caption{Maximal input for different values of~$c$ in~\eqref{def:vp}.}
\label{Fig:u_c}
\end{figure}
Figure~\ref{Fig:u_c} shows the relation between~$c$ and the maximal control input~$\|u\|_\infty$.
The maximal applied control increases with increasing~$c$.
This can be understood via the following reasoning.
First, in virtue of Remark~\ref{Rem:Comments-on-parameters}~iv), larger values of~$c$ cause larger lower bounds of the bijection~$\alpha(\cdot)$.
Second, the funnel boundary is 
$1/\vp(t) = c (T-t)$. So larger values of~$c$ result in a faster decay of the boundary,
and so the distance of the error to the funnel boundary decreases faster. 
Hence, larger input values are required to push the error away from the boundary.

\section{Conclusion}
We proposed a feedback controller, which achieves exact tracking in predefined finite time, while the tracking error evolves within prescribed boundaries.
We rigorously proved boundedness of all signals, and that the error as well as all its relevant derivatives vanish at the predefined final time.

\begin{ack}
I am deeply indebted to Thomas Berger (University of Paderborn) for many fruitful discussions.
I am also grateful for the comments of the anonymous reviewers who drew my attention to some important aspects.
\end{ack}

\section*{Appendix: Proof of Theorem~\ref{Thm:Exact-tracking-in-finte-time}}
First, we state the followin result.
\begin{Lemma} \label{Lem:x_le_M}
If for~$x \in \cC^1([\tau,T) ; \R^m)$, $m \in \N$, there exists~$M \ge \|x(\tau)\| \ge 0$ such that 
\begin{equation} \label{eq:ddt_x_le_0}
\forall \, t \in [\tau,T) \, : \ \left(
 \|x(t)\| \ge M \ \Rightarrow \ \ddt \|x(t)\|^2 \le 0 \right),
\end{equation}
then
\begin{equation} \label{eq:x_le_M}
\forall\, \ t \in [\tau,T) \, : \ \|x(t)\| \le M.
\end{equation}
\end{Lemma}
\begin{proof}
The proof follows the ideas in~\cite[Thm.~4.3]{Lanz21}.
Seeking a contradiction, we assume that there exists $t_1 \in (\tau,T)$ such that $\|x(t_1)\| > M$.
Then, by continuity, there exists $t_0 := \max \setdef{ t \in [\tau,t_1)}{ \|x(t)\| = M}$, and hence we have $\|x(t)\| \ge M$ for all $t \in [t_0,t_1]$.
Using \eqref{eq:ddt_x_le_0} we obtain
$
\|x(t_1)\|^2 - \|x(t_0)\|^2 = \int_{t_0}^{t_1} \left( \ddt \tfrac{1}{2} \|x(t)\|^2 \right) \dt \le 0,
$
and hence
$
M^2 < \|x(t_1)\|^2 \le \|x(t_0)\|^2 = M^2,
$
a contradiction.
Therefore,~\eqref{eq:x_le_M} holds for all~$t \in [\tau,T)$.
\end{proof}

\begin{proof}[Proof of Theorem~\ref{Thm:Exact-tracking-in-finte-time}]
The proof consists of eight steps. \\
\textit{Step one.} 
We show existence of a solution of~\eqref{eq:system},~\eqref{def:u}.
To this end, we aim to reformulate~\eqref{eq:system},~\eqref{def:u} as an initial value problem of the form
\begin{equation} \label{eq:IVP}
\begin{aligned}
\dot x(t) &= F\left(t,x(t), \textbf{T}(x)(t) \right), \\
x (0) &= \left(y^0(0), \dot y^0(0),\ldots, (\ddt)^{r-1}y^0(0) \right),
\end{aligned}
\end{equation}
where we set
$
x(\cdot) = ( y(\cdot), \dot y(\cdot), \ldots, y^{(r-1)}(\cdot) ),
$
and ${n=rm}$.
Setting~$\cD_0 := \setdef{ v \in \R^m}{\|v\| < 1}$ we choose some interval~$I \subseteq [0,T)$ with~$0 \in I$ such that $(e_1,\ldots,e_r) : I \to \R^{n}$ satisfy the relations in~\eqref{def:control_variables} and be such that for all~$t \in I$ we have $e_1(t),\ldots,e_{r-1}(t) \in \cD_0$, which is possible via the initial conditions~\eqref{eq:initial}.
Then, with the aid of~\eqref{def:Gam_k} we have for all $k=1,\ldots,r-1$
\begin{equation*}
\gamma_k^{(j)}(t) = \Gamma_j \big( t,e_k(t),\ldots,e_{k+j}(t) \big), \ 0 \le j \le r-k, \ t \in I.
\end{equation*}
Next, we define the function
\begin{equation*}
\begin{aligned}
\tilde e_1 : [0,T) \times \R^m &\to \R^m, \quad
(t,\xi_0)  \mapsto \vp(t) \left( \xi_0 - y_{\rm ref}(t) \right), 
\end{aligned}
\end{equation*}
and the set
$
\tilde \cD_1 := \setdef{ (t,\xi_0) \in [0,T) \times \R^m }{  \tilde e_1(t,\xi_0) \in \cD_0  } .
$
With this, we recursively define for~$k=2,\ldots,r$ the functions
\begin{equation*}
\begin{aligned}
\tilde e_k \, : \tilde \cD_{k-1} \times \R^m &\to \R^m, \\
(t,\xi_0,\ldots,\xi_{k-1}) & \mapsto \vp(t) \left( \xi_{k-1} - y_{\rm ref}^{(k-1)}(t) \right) \\
& \quad + \vp(t) \sum_{i=1}^{k-1} \Gamma_{k-1-i} \big( t,\tilde e_i,\ldots,\tilde e_{k-1} \big),
\end{aligned}
\end{equation*}
where for the sake of better legibility we omit the arguments of~$\tilde e_j$, $j=1,\ldots,k-1$; further, we define the sets
\begin{equation*}
\begin{aligned}
\tilde \cD_{k} := \big\{ &(t,\xi_0,\ldots,\xi_{k-1}) \in \tilde \cD_{k-1} \times \R^m  \ \vline \\
& \quad  \tilde e_{i}(t,\xi_0,\ldots,\xi_{k-1})  \in \cD_0, i=1,\ldots,k \big\} .
\end{aligned}
\end{equation*}
Together with~$\alpha(\cdot), N(\cdot)$ as in~\eqref{def:control_functions} 
we define for~$t \in I$
\begin{equation*}
N_r(t) := \big( N \circ \alpha \big) \left( \| \tilde e_r \big( t,y(t),\dot y(t),\ldots,y^{(r-1)}(t) \big) \|^2 \right).
\end{equation*}
Then, the control~$u(\cdot)$ defined in~\eqref{def:u} reads
\begin{equation*}
\begin{aligned}
u(t) = N_r(t) \cdot \tilde e_r \big( t,y(t),\ldots,y^{(r-1)}(t) \big), \ t \in I.
%
%
\end{aligned}
\end{equation*}
Lastly, we define the function~$F : \tilde \cD_{r-1} \times \R^q \to \R^{n}$ by
\begin{equation*}
\begin{aligned}
 (t,\xi_0,\ldots, \xi_{r-1}, \eta) \mapsto  \left( \xi_1,\ldots,\xi_{r-1}, f \left( d(t), \eta, N_r(t) \cdot \tilde e_r \right) \right),
\end{aligned}
\end{equation*}
where in $\tilde e_r = \tilde e_r(t,\xi_0,\ldots,\xi_{r-1})$ we omit the arguments for the sake of better legibility.
Together, the initial value problem~\eqref{eq:system},~\eqref{def:u} is equivalent to~\eqref{eq:IVP}.
In particular, we have~$(0,x(0)) \in \tilde \cD_{r-1}$, the function~$F$ is measurable in the variable~$t$, continuous in~$(\xi_0,\ldots,\xi_{r-1}, \eta)$ and locally essentially bounded.
Hence,~\cite[Thm.~B.1]{IlchRyan09} yields the existence of a maximal solution~$x: [-\sigma,\omega) \to \R^{n}$ of~\eqref{eq:IVP}, $0 < \omega \le T$.
In particular, the graph of the solution of~\eqref{eq:IVP} is not a compact subset of~$\tilde \cD_{r-1} $. 
\\
\textit{Step two.}
For the functions~$e_k(\cdot)$ introduced in~\eqref{def:e_k} we show 
that for all~$k=1,\ldots,r-1$ there exists~$\ve_k \in (0,1)$ such that~$\| e_k(t)\| \le \ve_k$ for all~$t \in [0,\omega)$.
We observe that for~$t \in [0,\omega)$ and~$k = 1,\ldots,r$ we have
\begin{equation*}
e_k(t) - \vp(t) \sum_{i=1}^{k-1} \gamma_{i}^{(k-1-i)}(t) = \vp(t) e^{(k-1)}(t).
\end{equation*}
We define $\alpha_k(t) := \alpha(\|e_k(t)\|^2)$.
Using~$ \dot \vp(t) = c \vp(t)^2$, we calculate for $k=1,\ldots,r-1$
\begin{equation} \label{eq:ddt_ek}
\begin{scriptsize}
\begin{aligned}
\dot e_k(t) &= \dot \vp(t) e^{(k-1)}(t) + \vp(t) e^{(k)}(t) \\
& + \dot \vp(t) \sum_{i=1}^{k-1} \gamma_{i}^{(k-1-i)}(t) + \vp(t) \sum_{i=1}^{k-1} \gamma_{i}^{(k-i)}(t) \\
& = \frac{\dot \vp(t)}{\vp(t)} \Big( e_k(t) - \vp(t) \sum_{i=1}^{k-1} \gamma_{i}^{(k-1-i)}(t) \Big) \\
&+ \Big( e_{k+1}(t) - \vp(t) \sum_{i=1}^{k} \gamma_{i}^{(k-i)}(t) \Big) \\
& + \dot \vp(t) \sum_{i=1}^{k-1} \gamma_{i}^{(k-1-i)}(t) + \vp(t) \sum_{i=1}^{k-1} \gamma_{i}^{(k-i)}(t) \\
& = (c-\alpha_k(t) ) \vp(t) e_k(t) + e_{k+1}(t),
%
%
\\
\dot e_r(t) &= c \vp(t) e_r(t) + \vp(t) e^{(r)}(t) +  \vp(t) \sum_{i=1}^{r-1} \gamma_{i}^{(r-i)}(t).
\end{aligned}
\end{scriptsize}
\end{equation}
Using the definitions of~$\alpha_k(\cdot)$ and~$\gamma_k(\cdot)$, we have
\begin{align} \label{eq:ddt_alp}
\dot \gamma_k(t) &= \ddt ( \alpha_k(t) e_k(t) ) \\
& = 2 \alpha'(\|e_k(t)\|^2) \langle e_k(t), \dot e_k(t) \rangle e_k(t) + \alpha_k(t) \dot e_k(t) . \nonumber
\end{align}
We observe~$e_k(t) = \tilde e_k(y(t),\ldots,y^{(k-1)}(t))$ and hence, since~$\tilde e_k(\cdot) \in \cD_0$ due to the initial conditions~\eqref{eq:initial}, we have
\begin{equation*}
\forall \, k=1,\ldots,r \ \forall \, t \in [0,\omega) \,: \ \|e_k(t)\| < 1.
\end{equation*}
We set $\hat \ve_k := \| e_k(0)\|^2 < 1$ and $\lambda := \vp(0) = \inf_{s \in [0,T)} \vp(s)> 0$.
Let~$\ve$ be the unique point in~$(0,1)$ such that $\alpha(\ve) \ve = (1 + c \lambda)/\lambda$ and define $\ve_k := \max\{\ve,\hat \ve_k\} < 1$.
We show that
\begin{equation} \label{eq:ek_le_vek}
\forall\, k\in\{1,\ldots,r-1\} \ \forall \, t \in [0,\omega) \, : \ \|e_k(t)\|^2 \le \ve_k.
\end{equation}
Seeking a contradiction, we suppose this is false for at least one~$\ell \in \{1,\ldots,r-1\}$. 
Then there exists $t_1 \in (0,\omega)$ such that 
$\|e_\ell(t_1)\|^2 > \ve_\ell$.
We define
$
t_0 := \max \setdef{ t \in [0,t_1)}{ \|e_\ell(t)\|^2 = \ve_\ell}.
$
Then we have
\begin{equation*}
\forall\, t \in [t_0,t_1]\,: \ \ve \le \ve_\ell \le \|e_\ell(t)\|^2,
\end{equation*}
which gives, invoking monotonicity of~$\alpha(\cdot)$,
that $\alpha(\ve) \le \alpha(\|e_\ell(t)\|^2) = \alpha_\ell(t)$
for all~$t \in  [t_0,t_1]$.
Hence, we have
$\alpha_\ell(t) \|e_\ell(t)\|^2 \ge \alpha(\ve) \ve = \frac{1+ c \lambda}{\lambda}$ 
for all~$t \in [t_0,t_1]$.
With this, using~$\alpha_\ell(\cdot) \ge c$ via~\eqref{def:alp}, and the relations in~\eqref{eq:ddt_ek}, we calculate for~$t \in [t_0,t_1]$
\begin{equation*}
\begin{aligned}
\ddt \tfrac{1}{2} \|e_\ell(t)\|^2 &= \langle e_\ell(t), (c-\alpha_\ell(t) ) \vp(t) e_\ell(t) + e_{\ell+1}(t) \rangle \\
& = \! - \vp(t) (\alpha_\ell(t) \!-\! c) \|e_\ell(t)\|^2 \! + \! \langle e_\ell(t), e_{\ell+1}(t) \rangle \\
& < -\lambda \alpha_\ell(t) \|e_{\ell}(t)\|^2 + c\lambda + 1 \le 0 ,
\end{aligned}
\end{equation*}
which implies
$
\ve_\ell < \|e_\ell(t_1)\|^2 < \|e_\ell(t_0)\|^2 = \ve_\ell,
$
a contradiction.
Therefore~\eqref{eq:ek_le_vek} holds.
This implies boundedness of~$\alpha_k$ (bounded by $\alpha(\ve_k)$) and boundedness of~$\gamma_k$ (bounded by $\alpha(\ve_k) \sqrt{\ve_k}$) for all $k=1,\ldots,r-1$. 
\\
\textit{Step three.}
Since the functions~$e_k(\cdot)$ defined in~\eqref{def:e_k} involve higher order derivatives of the functions~$\gamma_k(\cdot)$ we aim to show boundedness of the latter. 
To this end, recalling the definition of~$\gamma_k(\cdot)$, we establish boundedness of higher
order derivatives of~$\alpha_k(\cdot)$ on~$[0,\omega)$, which in turn involve higher order derivatives of~$e_k(\cdot)$.
Hence, we show boundedness of higher order derivatives of~$e_k(\cdot)$ on~$[0,\omega)$; more precise, we show boundedness of~$e_k^{(r-k)}(\cdot)$ on~$[0,\omega)$ for $k=1,\ldots,r-1$.
Recalling the definition of~$\vp(\cdot)$ in~\eqref{def:vp} we have $\vp^{(j)}(t) = c^{j} j! \vp(t)^{j+1}$ for~$j\in\N$. 
Using the generalized Leibniz rule, we obtain via~\eqref{eq:ddt_ek} for~$k=1,\ldots,r-1$ and~$1 \le j \le r-k$ the recursion
\begin{equation} \label{eq:j_th-derivative-ek}
\begin{scriptsize}
\begin{aligned}
e_k^{(j)}(t) &= \big( (c-\alpha_k(t)) \vp(t) e_k(t) \big)^{(j-1)} + e_{k+1}^{(j-1)}(t) \\
&= \sum_{j_1+j_2+j_3 = j-1} \frac{(j-1)!}{j_1! j_2! j_3!}  (c-\alpha_k(t))^{(j_1)} \vp^{(j_2)}(t) e_k^{(j_3)}(t) \\
& + e_{k+1}^{(j-1)}(t) \\
&= \sum_{j_1+j_2+j_3 = j-1} \frac{(j-1)!}{j_1! j_3!}  (c-\alpha_k(t))^{(j_1)} c^{j_2} \vp(t)^{j_2+1} e_k^{(j_3)}(t) \\
& + e_{k+1}^{(j-1)}(t) .
\end{aligned}
\end{scriptsize}
\end{equation}
We expatiate on the expression above for~$j=2$:
\begin{equation*}
\begin{aligned}
\ddot e_k(t) 
& =  -2\alpha'(\|e_k(t)\|^2) \langle e_k(t), \dot e_k(t) \rangle \vp(t) e_k(t) \\
& + (c-\alpha_k(t))c \vp^2 e_k(t)   \\
& + (c\!-\! \alpha_k(t)) \left( (c \!- \! \alpha_k(t)) \vp(t)^2 e_k(t) \!+\! \vp(t) e_{k+1}(t) \right) \\
& + (c-\alpha_{k+1}(t)) \vp(t) e_{k+1}(t) + e_{k+2}(t).
\end{aligned}
\end{equation*}
This recursion successively leads to the following observations. Since~$j_1+j_2+j_3 = j-1$ 
\begin{itemize}[leftmargin=1.4em]
\item for~$j_1=0$ the expression~$\vp(\cdot)^{j_2+1} e_k^{(j_3)}(\cdot)$ involves at most the~$j-1^{\rm st}$ derivative of~$e_k(\cdot)$, 
and at most the~$j^{\rm th}$ power of~$\vp(\cdot)$; 
the other terms involve (at most) derivatives and powers of the form~$\vp(\cdot)^{j_2+1} e_k^{(j-1-j_2)}(\cdot)$ for~$j\le r-k$,
\item $e_k^{(j)}(\cdot)$ involves $e_{k+1}^{(j-1)}(\cdot)$ which itself involves~$e_{k+2}^{(j-2)}(\cdot)$, therefore~$e_k^{(j)}(\cdot)$ involves~$e_{k+j}(\cdot)$,
\item the highest derivative of~$\alpha_k(\cdot)$ appearing in $e_k^{(j)}(\cdot)$ is~$\alpha_k^{(j-1)}(\cdot)$, which itself involves at most the~~$j-1^{\rm st}$ derivative of~$e_k(\cdot)$.
\end{itemize}
These observations together with the fact that
\begin{equation} \label{eq:vp+j_ek-vp+j-1_ek}
\begin{aligned}
\forall\, j \in \N\,:\ & \vp(\cdot)^{j} e_k(\cdot) \in \cL^{\infty}([0,\omega) ; \R^m) \\ 
& \Rightarrow   \vp(\cdot)^{j-1} e_k(\cdot) \in \cL^{\infty}([0,\omega) ; \R^m)
\end{aligned}
\end{equation}
yield that boundedness of~$e_k^{(j)}(\cdot)$ on~$[0,\omega)$ can be established by showing boundedness of~$\vp(\cdot)^{r-k} e_k(\cdot)$ for all~$k=1,\ldots,r-1$.
In order to show this, we initially establish the following:
for all~$k=2,\ldots,r-1$
\begin{equation} \label{eq:vp+r-k_e_k-bounded-vp+r-k+1_e_k-1-bounded}
\begin{aligned}
 & \vp(\cdot)^{r-k} e_k(\cdot) \in \cL^{\infty}([0,\omega) ; \R^m) \\ 
& \Rightarrow   \vp(\cdot)^{r-k+1} e_{k-1}(\cdot) \in \cL^{\infty}([0,\omega) ; \R^m).
\end{aligned}
\end{equation}
To see this, let $\vp(\cdot)^{r-k} e_k(\cdot) \in \cL^{\infty}([0,\omega) ; \R^m)$ and set $M_k := \sup_{s \in [0,\omega)} \| \vp(s)^{r-k} e_k(s)\| < \infty$. 
Then we consider for~$t \in [0,\omega)$
\begin{equation*}
\begin{scriptsize}
\begin{aligned}
& \ddt \tfrac{1}{2} \| \vp(t)^{r-k+1} e_{k-1}(t)\|^2  \\
&= \langle \vp(t)^{r-k+1} e_{k-1}(t), \vp(t)^{r-k+1} \dot e_{k-1}(t) \rangle \\
& \qquad + \langle \vp(t)^{r-k+1} e_{k-1}(t),  c(r-k+1) \vp(t)^{r-k+2} e_{k-1}(t) \rangle \\
& = \langle \vp(t)^{r-k+1} e_{k-1}(t), \vp(t)^{r-k+1} \big( (c-\alpha_{k-1}(t))\vp(t)e_{k-1}(t) \rangle \\
& \qquad + \langle \vp(t)^{r-k+1} e_{k-1}(t),  c(r-k+1) \vp(t)^{r-k+2} e_{k-1}(t)\rangle \\
& \qquad + \langle \vp(t)^{r-k+1} e_{k-1}(t), \vp(t)^{r-k+1} e_k(t) \rangle \\
& \le - \vp(t) \big( c(r+1) - c(r-k+2) \big) \|\vp(t)^{r-k+1} e_{k-1}(t)\|^2 \\
& \qquad + \vp(t) M_{k} \|\vp(t)^{r-k+1} e_{k-1}(t)\| \\
&\le -\vp(t) \Big( c(k-1) \|\vp(t)^{r-k+1} e_{k-1}(t)\| \\
& \qquad \qquad \qquad \qquad \qquad  - M_{k} \Big) \cdot  \|\vp(t)^{r-k+1} e_{k-1}(t)\| 
\end{aligned}
\end{scriptsize}
\end{equation*}
which is non-positive for~$\|\vp(t)^{r-k+1} e_{k-1}(t)\| \ge \tfrac{M_k}{c(k-1)}$. 
Hence, Lemma~\ref{Lem:x_le_M} yields boundedness of~$\vp(\cdot)^{r-k+1} e_{k-1}(\cdot)$ on~$[0,\omega)$. 
A successive application of~\eqref{eq:vp+r-k_e_k-bounded-vp+r-k+1_e_k-1-bounded} yields
\begin{equation} \label{eq:vp+r-k-e_k-bounded}
\!\! \forall \, k=1,\ldots,r-1\, : \ \vp(\cdot)^{r-k} e_k(\cdot) \in \cL^\infty([0,\omega) ; \R^m).
\end{equation}
In particular, via~\eqref{eq:vp+j_ek-vp+j-1_ek} we have~$\vp(\cdot) e_{k}(\cdot) \in \cL^\infty([0,\omega) ; \R^m)$ for all~$k=1,\ldots,r-1$ from which boundedness of~$\dot e_k(\cdot)$ follows, which in turn implies boundedness of~$\dot \alpha_k(\cdot)$ and~$\dot \gamma_k(\cdot)$ on~$[0,\omega)$.
By~\eqref{eq:j_th-derivative-ek} and~\eqref{eq:vp+j_ek-vp+j-1_ek} boundedness of~$e^{(j)}_k(\cdot)$ successively follows for all~$j \le r-k-1$, from which we may deduce boundedness of~$ \alpha_k^{(j)}(\cdot)$ and~$ \gamma_k^{(j)}(\cdot)$ for~$j \le r-k-1$. 
Thus, for all~$k=1\ldots,r-1$ and~$j\le r-k-1$ there exists $ \bar \gamma_k^{j} := \sup_{s \in [0,\omega)} \gamma_k^{(j)}(s) < \infty $.
\\
\textit{Step four.} We show boundedness of~$x(\cdot)$ on~$[0,\omega)$.
Recalling the definition of~$e_k(\cdot)$ we see that for all $k=1,\ldots,r$ we have via the previous steps
\begin{equation} \label{eq:e+k-1-bounded}
    \begin{scriptsize}
\begin{aligned} 
\! \! \! \forall\, t \in [0,\omega)  : \, \|e^{(k-1)}(t)\| & \le \left\|\frac{e_{k}(t)}{\vp(t)} \right\| + \left\|\sum_{i=1}^{k-1} \gamma_i^{(k-1-i)}(t) \right\|  \\
 & \le  \frac{1}{\lambda}  + \sum_{i=1}^{k-1} \bar \gamma_i^{k-1-i} < \infty. 
\end{aligned}
\end{scriptsize}
\end{equation}
Therefore, since $x(\cdot) = (y(\cdot),\ldots,y^{(r-1)}(\cdot)) = (e(\cdot) + y_{\rm ref}(\cdot), \ldots, e^{(r-1)}(\cdot) + y_{\rm ref}^{(r-1)}(\cdot)) $ 
and the reference $y_{\rm ref} \in \cW^{r,\infty}([0,T) ; \R^m)$, we have $x \in \cL^\infty([0,\omega) ; \R^{rm})$.
\\
\textit{Step five.}
We show boundedness of~$\alpha_r(\cdot)$ on~$[0,\omega)$. 
Invoking the previous steps, in particular boundedness of~$x(\cdot)$ on~$[0,\omega)$, and the properties of the operator class $\cT^{n,q}_\sigma$ we deduce the existence of a compact $K_q \subset \R^q$ such that $T(x)(t) \in K_q$ for~$t \in [0,\omega)$; furthermore, since $d \in L^\infty(\rp ; \R^p)$ there exists a compact $K_p \subset \R^p$ such that $d(t) \in K_p$ for~$t \in [0,\omega)$.
By the high-gain property there exists~$\rho \in (0,1)$ such that the continuous function
$\chi(s) = \min \big\{ \langle v,f(\delta,z,-sv) \rangle  \, \vline \, (\delta,z,v) \in K_q \!\times\! K_p \!\times\! V \big\}$
is unbounded from above; 
where we define the compact set
$
V := \setdef{ v \in \R^m}{ \rho \le \|v\| \le 1}.
$
We show boundedness of~$\alpha_r(\cdot)$ by contradiction.
Since
$N : \rp \to \R$ is surjective, the set~$\setdef{ \kappa > \rho_0 }{ N(\kappa) = \rho_1}$ is non-empty for every~$\rho_0 \in \rp$ and every~$\rho_1 \in \R$.
Following the proof in~\cite[pp.~188-190]{BergIlch21},
we choose a sequence~$(s_j)$ such that the corresponding sequence~$\left(\chi(s_j)\right)$ is positive, strictly increasing and unbounded.
We initialize a sequence~$(\kappa_j)$ by $\kappa_1 > \alpha(\rho^2) + \alpha_r(0)$ such that~$N(\kappa_1) = s_1$, and hereinafter define the strictly increasing sequence~$(\kappa_j)$ by
$
 \kappa_{j+1} := \inf \setdef{ \kappa > \kappa_j }{ N(\kappa) = s_{j+1} },
$
which yields that
$
\lim_{j \to \infty} \chi(N(\kappa_j)) = \lim_{j \to \infty} \chi(s_j) = \infty.
$
Now, since we assumed~$\alpha_r(\cdot)$ to be unbounded and we have~$\kappa_{j+1} > \kappa_1 > \alpha_r(0)$ for all~$j \in \N$, we may define the sequence
$
\tau_j := \inf \setdef{ t \in [0,\omega) }{ \alpha_r(t) = \kappa_{j+1}  }, \ j \in \N_0,
$
which lies within~$(0,\omega)$.
Note that~$(\tau_j)$ is strictly increasing and we have $N(\alpha_r(\tau_j)) = N(\kappa_{j+1}) = s_{j+1}$ for each~$j \in \N_0$.
Next, we define a second sequence in~$(0,\omega)$
$
\sigma_j = \sup \setdef{ t \in [\tau_{j-1} , \tau_j]  }{ \chi(N(\alpha_r(t))) = \chi(s_j)  }, \ j \in \N.
$
%
With this, since the sequence~$(\chi(s_j))$ is strictly increasing, we obtain
$
 \chi(N(\alpha_r(\sigma_j))) = \chi(s_j) < \chi(s_{j+1}) = \chi(N(\alpha_r(\tau_j)))
$
for all~$j \in \N$,
and therefore, for all $j \in \N$ we have
$
\sigma_j < \tau_j,
$
and for all $t \in (\sigma_j, \tau_j]$ we have $\chi(N(\alpha_r(\sigma_j))) = \chi(s_j)    <  \chi(N(\alpha_r(t)))$.
Next, we show by contradiction that for all~$j \in \N$ and for all~$t \in [\sigma_j,\tau_j]$ we have~$e_r(t) \in V$. 
To this end, we first show -  by contradiction - that for all~$j \in \N$ we have~$\alpha_r(t) \ge \kappa_j$ for~$t \in [\sigma_j, \tau_j]$.
Suppose~$\alpha_r(t) < \kappa_j$ for some~$t \in [\sigma_j,\tau_j]$. 
Then by~$\alpha_r(\tau_j) = \kappa_{j+1} > \kappa_j$ and by continuity of~$\alpha_r$ there exists~$\tilde t \in (\sigma_j, \tau_j)$ such that~$\alpha_r(\tilde t) = \kappa_j$.
Hence, we find
$
\chi(N(\alpha_r(\tilde t))) = \chi(N(\kappa_j)) = \chi(s_j),
$
which contradicts the definition of~$\sigma_j$.
Therefore, $\alpha_r(t) \ge \kappa_j$ for all $t \in [\sigma_j,\tau_j]$.
Suppose~$e_r(t) \notin V$, i.e., since for all $t \in [0,\omega)$ we have $\|e_r(t)\| < 1$, suppose $\|e_r(t)\| < \rho$ for some $[\sigma_j,\tau_j]$.
This, together with $\alpha_r(t) \ge \kappa_j$, leads to the contradiction
$
\alpha(\rho^2) < \kappa_1 \le \kappa_j \le \alpha_r(t) = \alpha(\|e_r(t)\|^2) < \alpha(\rho^2).
$
Hence, we deduce $e_r(t) \in V$ for all $t \in [\sigma_j, \tau_j]$ and all $j \in \N$.
Since~$d(t) \in K_p$ and~$T(x)(t) \in K_q$ for~$t \in [0,\omega)$ we obtain, using $\sigma_j < \tau_j$, for all~$j \in \N$ and~$t \in [\sigma_j, \tau_j]$
\begin{equation}
\begin{scriptsize}    
\begin{aligned} \label{eq:er_pari_f}
&\langle e_r(t),   f\big( d(t), (Tx)(t), u(t) \big) \rangle \nonumber \\
& = - \langle -e_r(t) , f\big( d(t), \textbf{T}(x)(t), -N(\alpha_r(t)) (-e_r(t)) \rangle \nonumber \\
& \le - \min \setdef{ \langle v, f(\delta, z, -N(\alpha_r(t)) v }{ \begin{array}{l}
(\delta , z, v) \\ \in K_p \times K_q \times V
\end{array}  } \nonumber \\
& = - \chi(N(\alpha_r(t))) 
 \le - \chi(s_j) .
\end{aligned}
\end{scriptsize}
\end{equation}
Since~$y_{\rm ref} \in \cW^{r,\infty}([0,T) ; \R^m)$, we may set $c_{\rm ref} := \sup_{s \ge 0}\|y_{\rm ref}^{(r)}(s)\| < \infty$, and we recall $\sum_{i=1}^{r-1} \bar \gamma_i^{r-i} < \infty$ from the previous steps.
Furthermore, we observe~$\sigma_1 > 0$ and thus, by properties of~$\vp(\cdot)$, we may define $0 < \inf_{s \in [\sigma_1, T)} \vp(s) =: c_{\vp}$.
Then, with the aid of~\eqref{eq:ddt_ek} and~\eqref{eq:er_pari_f} for all~$j \in \N$ and~$t \in [\sigma_j, \tau_j]$ we obtain
\begin{equation*}
\begin{scriptsize}
\begin{aligned}
\ddt \tfrac{1}{2} \|e_r(t)\|^2 &= \langle e_r(t), c \vp(t) e_r(t) \rangle \\
&    + \langle e_r(t) \vp(t) \big( f(d(t), (Tx(t), u(t)) \!-\! y_{\rm ref}^{(r)}(t) \big) \rangle \\
&  + \langle e_r(t), \vp(t) \sum_{i=1}^{r-1} \gamma_i^{(r-i)}(t) \rangle \\
& < \vp(t) \Big( c + c_{\rm ref} + \sum_{i=1}^{r-1} \bar \gamma_i^{r-i} - \chi(s_j) \Big).
\end{aligned}
\end{scriptsize}
\vspace*{-1em}
\end{equation*}
Thus, still seeking a contradiction, we may choose~$J \in \N$ large enough such that for~$t \in [\sigma_J, \tau_J]$ we have
\begin{equation*}
\begin{scriptsize}
\begin{aligned}
& \vp(t) \Big( c + c_{\rm ref} + \sum_{i=1}^{r-1} \bar \gamma_i^{r-i} - \chi(s_J) \Big) \\
\vspace*{-2em}
& \le - c_{\vp} \Big(  \chi(s_J) - \Big(c + c_{\rm ref} + \sum_{i=1}^{r-1} \bar \gamma_i^{r-i} \Big)\Big) < 0,
\end{aligned}
\end{scriptsize}
\vspace*{-1em}
\end{equation*}
which yields
$
\|e_r(\tau_J)\|^2 < \| e_r(\sigma_J)\|^2,
$
which in turn gives 
$\alpha_r(\tau_J) = \alpha(\|e_r(\tau_J)\|^2) < \alpha(\|e_r(\sigma_J)\|^2) = \alpha_r(\sigma_J)$
for $t \in [\sigma_J, \tau_J]$.
This contradicts the definition of~$\tau_J$, by which we have~$\alpha_r(t) < \alpha_r(\tau_J)$ for all~$t \in [0,\tau_J)$.
Hence, the assumption of an unbounded~$\alpha_r(\cdot)$ cannot be true. 
As a direct consequence thereof, we may infer the existence of~$\ve_r \in (0,1)$, such that
$\|e_r(t)\| \le \ve_r$ for all $t \in [0,\omega)$.
\\
\textit{Step six.} We show~$\omega = T$. 
Via the previous steps we have for all $k=1,\ldots,r$ and all $t \in [0,\omega)$ that $\|e_k(t)\| \le  \ve^* := \sqrt{\max\{\ve_1,\ldots,\ve_r\}} < 1$,
by which the set
$
\hat \cD :=  \setdef{ (\zeta_1,\ldots,\zeta_r) \in \R^{rm} }{ \|\zeta_i\| \le \ve^*, i=1,\ldots,r} 
$
is a compact subset of~$\tilde \cD_{r-1} $.
Assume~$\omega < T$. Then, $x(t) \in \hat \cD \subset \tilde \cD_{r-1}$ for all $t \in [0,\omega)$.
By compactness of~$\hat \cD$, the closure of the graph of the solution~$x(\cdot)$ of~\eqref{eq:IVP} on $[0,\omega)$ is a compact subset of $\tilde \cD_{r-1} $, which contradicts the findings of \textit{Step one}. Thus, $\omega = T$.
\\
\textit{Step seven.} Assertion~$ii)$ is a direct consequence of \textit{Step four} and \textit{Step six}; and assertion~$iii)$ follows from \textit{Step two} and \textit{Step six}.
\\
\textit{Step eight.}
We show that the tracking error~$e(\cdot)$ and its derivatives tend to zero as~$t \to T$, this is, we show
\begin{equation} \label{eq:e-exact}
\forall\, k=1,\ldots,r \, : \ \lim_{t\to T} \|e^{(k-1)}(t)\| = 0.
\end{equation}
Estimate~\eqref{eq:e+k-1-bounded} is too rough to show~$\eqref{eq:e-exact}$. 
Recalling the definition $\gamma_k(\cdot) = \alpha_k(\cdot) e_k(\cdot)$, and exemplary its derivative~\eqref{eq:ddt_alp}, we see that by \textit{Step three} not only $\gamma_k^{(j)}(\cdot)$ is bounded on~$[0,\omega)$ for $j\le r-k-1$ but with the aid of~\eqref{eq:vp+r-k-e_k-bounded} even the product 
$\vp(\cdot) \, \gamma_k^{(j)}(\cdot)$ is bounded on~$[0,\omega)$, 
i.e., for all $\ell=1,\ldots,r-1$ there exists 
$ \hat \gamma_\ell^{j} := \sup_{s \in [0,\omega)} \vp(s) \gamma_\ell^{(j)}(s) < \infty $ for $0 \le j \le r-\ell - 1$.
Invoking \textit{Step two} and \textit{Step five}, we may improve~\eqref{eq:e+k-1-bounded} for all~$k=1,\ldots,r$ and all~$t \in [0, \omega) $ as follows
\begin{equation*}
\begin{scriptsize}
\begin{aligned}
\|e^{(k-1)}(t)\| &\le \frac{\|e_{k}(t)\|}{\vp(t)}  + \frac{1}{\vp(t)} \left\|\sum_{i=1}^{k-1} \vp(t)\gamma_i^{(k-1-i)}(t) \right\|  \\
& \le  \frac{\sqrt{\ve_k} + \sum_{i=1}^{k-1} \hat \gamma_i^{k-1-i} }{\vp(t)}. \\
\end{aligned}
\end{scriptsize}
\end{equation*}
Since~$\omega=T$ by \textit{Step six}, and~$\lim_{t\to T} \vp(t) = \infty$, we obtain~\eqref{eq:e-exact} for all~$k=1,\ldots,r$, which shows assertion~$iv)$ of the theorem
and completes the proof.
\end{proof}

\bibliography{MST,NewReferences}

\begin{thebibliography}{10}

\bibitem{Andr08}
Vincent Andrieu, Laurent Praly, and Alessandro Astolfi.
\newblock Homogeneous approximation, recursive observer design, and output
  feedback.
\newblock {\em {SIAM} Journal on Control and Optimization}, 47(4):1814--1850,
  2008.

\bibitem{BasiPana16}
Michael Basin, Chandrasekhara~Bharath Panathula, Yuri~B. Shtessel, and Pablo
  Rodriguez-Ramirez.
\newblock Continuous finite-time higher-order output regulators for systems
  with unmatched unbounded disturbances.
\newblock {\em {IEEE} Transactions on Industrial Electronics}, pages 1--1,
  2016.

\bibitem{BergIlch21}
Thomas Berger, Achim Ilchmann, and Eugene~P. Ryan.
\newblock Funnel control of nonlinear systems.
\newblock {\em Mathematics of Control, Signals, and Systems}, 33(1):151--194,
  feb 2021.

\bibitem{BhatBern05}
S.~P. Bhat and D.~S. Bernstein.
\newblock Geometric homogeneity with applications to finite-time stability.
\newblock {\em Mathematics of Control, Signals, and Systems}, 17(2):101--127,
  may 2005.

\bibitem{ByrnIsid91a}
Christopher~I. Byrnes and Alberto Isidori.
\newblock Asymptotic stabilization of minimum phase nonlinear systems.
\newblock {\em {IEEE} Trans. Autom. Control}, 36(10):1122--1137, 1991.

\bibitem{ByrnWill84}
Christopher~I. Byrnes and Jan~C. Willems.
\newblock Adaptive stabilization of multivariable linear systems.
\newblock In {\em Proc. 23rd~{IEEE} Conf. Decis. Control}, pages 1574--1577,
  1984.

\bibitem{ClohWilt60}
W.~H. Clohessy and R.~S. Wiltshire.
\newblock Terminal guidance system for satellite rendezvous.
\newblock {\em Journal of the Aerospace Sciences}, 27(9):653--658, sep 1960.

\bibitem{Davi13}
Jorge Davila.
\newblock Exact tracking using backstepping control design and high-order
  sliding modes.
\newblock {\em {IEEE} Transactions on Automatic Control}, 58(8):2077--2081, aug
  2013.

\bibitem{DimaBech20}
Ioannis~S. Dimanidis, Charalampos~P. Bechlioulis, and George~A. Rovithakis.
\newblock Output feedback approximation-free prescribed performance tracking
  control for uncertain {MIMO} nonlinear systems.
\newblock {\em {IEEE} Transactions on Automatic Control}, 65(12):5058--5069,
  dec 2020.

\bibitem{EspiPerr22}
Nicolas Espitia and Wilfrid Perruquetti.
\newblock Predictor-feedback prescribed-time stabilization of {LTI} systems
  with input delay.
\newblock {\em {IEEE} Transactions on Automatic Control}, 67(6):2784--2799,
  2022-06.

\bibitem{EstrFrid16}
Antonio Estrada, Leonid Fridman, and Rafael Iriarte.
\newblock Combined backstepping and {HOSM} control design for a class of
  nonlinear {MIMO} systems.
\newblock {\em International Journal of Robust and Nonlinear Control},
  27(4):566--581, jul 2016.

\bibitem{Hill78}
G.~W. Hill.
\newblock Researches in the lunar theory.
\newblock {\em American Journal of Mathematics}, 1(3):245, 1878.

\bibitem{IlchRyan09}
Achim Ilchmann and Eugene~P. Ryan.
\newblock Performance funnels and tracking control.
\newblock {\em Int. J. Control}, 82(10):1828--1840, 2009.

\bibitem{IlchRyan02b}
Achim Ilchmann, Eugene~P. Ryan, and Christopher~J. Sangwin.
\newblock Tracking with prescribed transient behaviour.
\newblock {\em ESAIM: Control, Optimisation and Calculus of Variations},
  7:471--493, 2002.

\bibitem{Isid95}
Alberto Isidori.
\newblock {\em Nonlinear Control Systems}.
\newblock Communications and Control Engineering Series. Springer-Verlag,
  Berlin, 3rd edition, 1995.

\bibitem{RodrSanc17}
Esteban Jim{\'{e}}nez-Rodr{\'{\i}}guez, Juan~Diego S{\'{a}}nchez-Torres, and
  Alexander~G. Loukianov.
\newblock Predefined-time backstepping control for tracking a class of
  mechanical systems.
\newblock {\em {IFAC}-{PapersOnLine}}, 50(1):1680--1685, jul 2017.

\bibitem{Klue99}
Craig~A. Kluever.
\newblock Feedback control for spacecraft rendezvous and docking.
\newblock {\em Journal of Guidance, Control, and Dynamics}, 22(4):609--611, jul
  1999.

\bibitem{Lanz21}
Lukas Lanza.
\newblock Internal dynamics of multibody systems.
\newblock {\em Systems {\&} Control Letters}, 152:104931, jun 2021.

\bibitem{LeeTren19}
Jin~Gyu Lee and Stephan Trenn.
\newblock Asymptotic tracking via funnel control.
\newblock In {\em 2019 {IEEE} 58th Conference on Decision and Control ({CDC})}.
  {IEEE}, dec 2019.

\bibitem{Leva03}
Arie Levant.
\newblock Higher-order sliding modes, differentiation and output-feedback
  control.
\newblock {\em International Journal of Control}, 76(9-10):924--941, 2003.

\bibitem{Leva14}
Arie Levant.
\newblock Finite-time stabilization of uncertain {MIMO} systems.
\newblock In {\em 53rd {IEEE} Conference on Decision and Control}. {IEEE},
  2014-12.

\bibitem{MercUrib21}
Angel Mercado-Uribe, Jaime~A. Moreno, Andrey Polyakov, and Denis Efimov.
\newblock Multiple-input multiple-output homogeneous integral control design
  using the implicit lyapunov function approach.
\newblock {\em International Journal of Robust and Nonlinear Control},
  31(9):3417--3438, 2021.

\bibitem{NunePeix14}
Eduardo~V.L. Nunes, Alessandro~J. Peixoto, Tiago~Roux Oliveira, and Liu Hsu.
\newblock Global exact tracking for uncertain {MIMO} linear systems by output
  feedback sliding mode control.
\newblock {\em Journal of the Franklin Institute}, 351(4):2015--2032, 2014-04.

\bibitem{NuneHsu09}
E.V.L. Nunes, Liu Hsu, and F.~Lizarralde.
\newblock Global exact tracking for uncertain systems using output-feedback
  sliding mode control.
\newblock {\em {IEEE} Transactions on Automatic Control}, 54(5):1141--1147,
  2009-05.

\bibitem{OlivPeix13}
Tiago~Roux Oliveira, Alessandro~Jacoud Peixoto, and Liu Hsu.
\newblock Peaking free output-feedback exact tracking of uncertain nonlinear
  systems via dwell-time and norm observers.
\newblock {\em Int. J. Robust \& Nonlinear Control}, 23(5):483--513, 2013.

\bibitem{RyanSang09}
Eugene~P. Ryan, Christopher~J. Sangwin, and Philip Townsend.
\newblock Controlled functional differential equations: approximate and exact
  asymptotic tracking with prescribed transient performance.
\newblock {\em ESAIM: Control, Optimisation and Calculus of Variations},
  15:745--762, 2009.

\bibitem{schwartz1999global}
Ben Schwartz, Alberto Isidori, and Tzyh~Jong Tarn.
\newblock Global normal forms for {MIMO no}nlinear systems, with applications
  to stabilization and disturbance attenuation.
\newblock {\em Mathematics of Control, Signals and Systems}, 12:121--142, 1999.

\bibitem{song2017time}
Yongduan Song, Yujuan Wang, John Holloway, and Miroslav Krstic.
\newblock Time-varying feedback for regulation of normal-form nonlinear systems
  in prescribed finite time.
\newblock {\em Automatica}, 83:243--251, 2017.

\bibitem{Tren19}
Stephan Trenn.
\newblock Asymptotic tracking with funnel control.
\newblock {\em {PAMM}}, 19(1), sep 2019.

\bibitem{VergDima19}
Christos~K. Verginis and Dimos~V. Dimarogonas.
\newblock Asymptotic stability of uncertain {Lagrangian} systems with
  prescribed transient response.
\newblock In {\em Proceeding 58th IEEE Conference on Decision and Control,
  Nice, France}, pages 7037--7042. {IEEE}, 2019-12.

\bibitem{VergDima21}
Christos~K. Verginis and Dimos~V. Dimarogonas.
\newblock Asymptotic tracking of second-order nonsmooth feedback stabilizable
  unknown systems with prescribed transient response.
\newblock {\em IEEE Transactions on Automatic Control}, 66(7):3296--3302,
  2021-07.

\bibitem{YangDing18}
Jun Yang, Zhengtao Ding, Shihua Li, and Chuanlin Zhang.
\newblock Continuous finite-time output regulation of nonlinear systems with
  unmatched time-varying disturbances.
\newblock {\em IEEE Control Systems Letters}, 2(1):97--102, 2018-01.

\end{thebibliography}

\begin{wrapfigure}{l}{65pt}
\centering
\vspace*{-0.9em}
  \includegraphics[height=90pt]{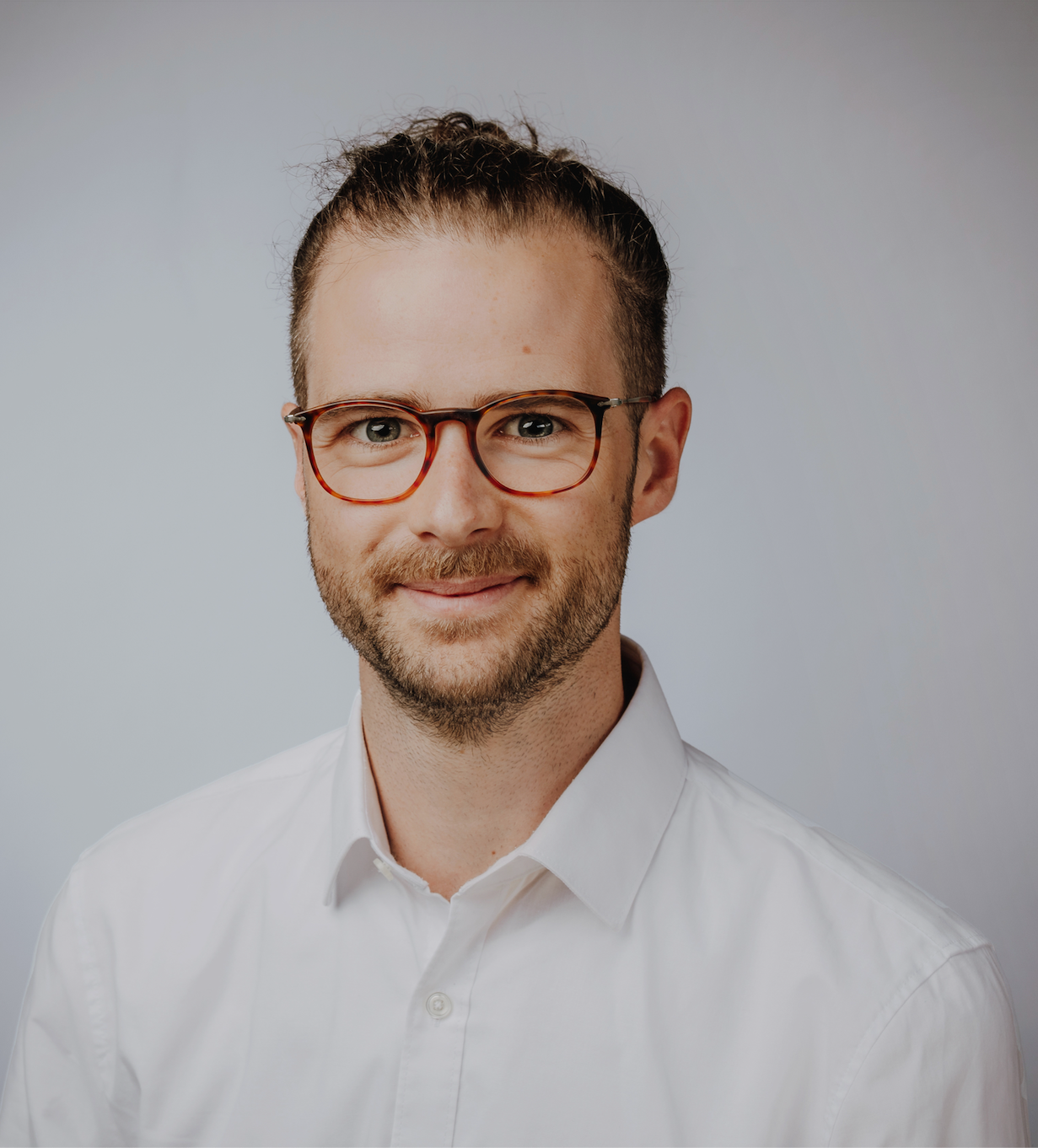} \par 
  \vspace*{-1em}
\end{wrapfigure}
\scriptsize
\noindent
\textbf{Lukas Lanza} received his Bachelor’s degrees in Physics and Philosophy from the University of M\"unster in 2015, his Master’s degree in Industrial Mathematics from the University of Hamburg in 2019, and his PhD in Mathematics from the University of Paderborn in 2022.
Since October 2022, he has been working as a Postdoc in the group Optimization-based Control at TU Ilmenau. Currently, his research focuses on safe data-driven control, predictive control with guarantees, and learning-based control, with applications to maglev trains, electrical drives, and energy systems.

\end{document}